\documentclass[a4paper,11pt]{article}
\usepackage{amsfonts}
\usepackage{mathrsfs}
\usepackage{amsmath,amssymb,amsthm}
\usepackage{latexsym}
\usepackage{indentfirst}
\usepackage{color}
\usepackage{footmisc}
\usepackage{endnotes}
\usepackage{algorithm}
\usepackage{float}
\usepackage{graphicx}
\usepackage{longtable}
\usepackage[authoryear,round]{natbib}
\usepackage[colorlinks,citecolor=blue,urlcolor=blue,filecolor=blue,backref=page]{hyperref}
\usepackage{graphicx}
\usepackage{tikz}
\usepackage[top=1in, bottom=1in, left=1in, right=1in]{geometry}

\DeclareMathOperator*{\E}{\mbox{E}}
\DeclareMathOperator*{\Jw}{\mathcal{J}_1}
\DeclareMathOperator*{\Jt}{\mathcal{J}_0}

\DeclareMathOperator*{\dett}{\text{det}}
\DeclareMathOperator*{\voll}{\text{Vol}}
\DeclareMathOperator*{\DD}{\bf D}

\newtheorem{thm}{Theorem}
\newtheorem{prop}{Proposition}
\newtheorem{property}{Property}

\begin{document}

\title{\textbf{On Bayesian Oracle Properties}}

\author{Wenxin Jiang \thanks{\noindent Shandong University (Taishan Scholar Adjunct Professor) and Northwestern University, \newline
\texttt{wjiang@northwestern.edu}}
\and
Cheng Li \thanks{\noindent National University of Singapore, \texttt{stalic@nus.edu.sg}}
}


\date{}
\maketitle

%
%
%
%
%



\begin{abstract}
When model uncertainty is handled by Bayesian model averaging (BMA) or Bayesian model selection (BMS), the posterior distribution possesses a desirable ``oracle property'' for parametric inference, if for large enough data it is nearly as good as the oracle posterior, obtained by assuming unrealistically that the true model is known and only the true model is used. We study the oracle properties in a very general context of quasi-posterior,  which can accommodate non-regular models with cubic root asymptotics and partial identification. Our approach for proving the oracle properties is based on a unified treatment that bounds the posterior probability of model mis-selection. This theoretical framework can be of interest to Bayesian statisticians who would like to theoretically justify their new model selection or model averaging methods in addition to empirical results. Furthermore, for non-regular models, we obtain nontrivial conclusions on the choice of prior penalty on model complexity, the temperature parameter of the quasi-posterior, and the advantage of BMA over BMS.
\end{abstract}


MSC2010 Classification Codes: 62E99, 62F15.

%


\section{Introduction}

The terminology of frequentist {\it oracle property} was first introduced in \citet{fl01} for a frequentist penalization method in model selection, by which statistical inferences ``work as well as if the correct submodel were known." Thereafter the oracle property has become a popular concept in the statistics literature. On the other hand, analogs of such an oracle property have not been widely studied in the Bayesian context, with the exception of a few recent works in special model setups (\citealt{ir11}, \citealt*{csv15}, \citealt{lj16}, etc.)

In this paper, we define different versions of Bayesian oracle properties in a general framework with quasi-posteriors and present a systematic way to study them by bounding the probability of model mis-selection. In particular, we are interested in the interplay between several different subjects: Bayesian model averaging (BMA), Bayesian model selection (BMS) based on the {\it Maximum-A-Posteriori} (MAP) model, and Bayesian posterior inference based on the unknown true model (i.e. the oracle model). We reveal some surprisingly simple and general relations between these different topics, and discuss their applications in non-regular models with cubic root asymptotics and partial identification.


We first introduce the basic notation we will use throughout this paper. Let $\DD$ be the observed data with sample size $n$. Let $M$ be a generic model index, and the “true model” $M^*$ be a possible value of $M$ which is related to the data generating mechanism. In Bayesian model averaging and model selection, we always consider a countable sequence of models $\{M_j\}$ indexed by $j=1,2,\ldots$, among which is the true model $M^*$. A prior probability $\pi(M_j)$ is assigned to each model $M_j$. Then each model $M_j$ proposes a different prior density $\pi(\theta |M_j)$ for the parameter $\theta$, supported on a parameter space $\Theta_j$, which can possibly overlap. The full parameter space is $\Theta=\cup_{j\geq 1}\Theta_j$. The overall prior distribution with density $\pi(\theta)$ is given by
\begin{align*}
&\pi(\theta) = \sum_{j\geq 1} \pi(\theta |M_j)\pi(M_j), \quad \text{for }\theta \in \Theta.
\end{align*}
Given the model $M_j$ and its proposed parameter $\theta$, let $p(\DD|\theta,M_j)$ be the likelihood function. Then the posterior density of $\theta$ through Bayesian model averaging (BMA) is given by
\begin{align*}
&\pi(\theta|\DD) \propto \sum_{j\geq 1} p(\DD|\theta,M_j)\pi(\theta |M_j)\pi(M_j), \quad \text{for }\theta \in \Theta.
\end{align*}
Throughout the paper, we use $\Pi$ to denote the underlying probability measure associated with density $\pi$.

Below we explain why the Bayesian version of oracle properties is desirable for dimension reduction in standard regular models, why the more general quasi-Bayesian framework is useful, and why our work will be of interest to the community of Bayesian statisticians.

\subsection{Bayesian oracle property is desirable for dimension reduction}

Consider a simple example of linear regression with known error variance, where $y\sim N\left(\sum_{j=1}^{p} x_j\theta_j,1\right)$, and $N(\mu,\sigma^2)$ denotes the normal distribution with mean $\mu$ and variance $\sigma^2$. Suppose that there exists an unknown true model $M^*$, in which only the first $p^*$ components of $\theta=(\theta_1,...,\theta_p)^\top$ are nonzero. Suppose we consider these nested candidate models $M_1,\ldots,M_p$, where the first $j$ components of $\theta$ are nonzero if $\theta$ comes from the model $M_j$. Given an observed independent and identically distributed (i.i.d.) sample $\DD=\left\{(y_i,x_{i1},...,x_{ip}),~ i=1,\ldots,n\right\}$, the BMA involves using the posterior
\begin{align*}
\pi\left(\theta |\DD \right)\propto \sum_{j=1}^{p} e^{-\frac{1}{2}\sum_{i=1}^n(y_i-\sum_{\ell=1}^j x_{i\ell}\theta_{\ell})^2} \pi\left(\theta|M_j\right) \pi\left(M_j\right).
\end{align*}
When $p\ll n$, we can set the prior $\pi(\theta|M_j)$ to be a component-wise independent product of normal priors, and $\pi(M_j)=1/p$ as a uniform prior.

For this simple example, the Bayesian oracle property can be roughly described as
$$\pi(\theta |\DD)\approx \pi(\theta|\DD,M^*) \propto e^{-\frac{1}{2}\sum_{i=1}^{n} (y_i-\sum_{\ell=1}^{p^*} x_{i\ell}\theta_{\ell})^2} \pi\left(\theta|M^*\right),$$
which is the posterior based on the true model $M^*$, as if we knew the truth $M^*$. This approximation can be in the sense of total variation norm, or in some other sense depending on what is regarded as meaningful.  This kind of result is desirable for automatic dimension reduction and variance reduction. If $p=10$ but $p^*=1$, then the mean squared error for  estimating the mean function $\sum_{j=1}^{p} x_j\theta_j$ can be reduced from about $10/n$ when using the full model $M_p$, to about $1/n$ when using BMA. Such advantage of BMA in dimension reduction and better prediction error has been empirically noticed in a variety of applications, such as in \citet{lj16} in the context of Bayesian generalized method of moments. When $p\gg n$, such dimension reduction through BMA is almost indispensable for any useful statistical inference, and has been widely studied in the literature with sparsity-inducing priors (\citealt{jr12}, \citealt{lsy13}, \citealt*{csv15}, etc.)

\subsection{It is useful to extend consideration to quasi-posteriors}

Our current paper extends the standard BMA to the general case of a {\it quasi-posterior}, where
$$\pi(\theta |\DD)\propto \sum_{j\geq 1} e^{-\lambda R_n(\DD|\theta, M_j)} \pi\left(\theta|M_j\right) \pi\left(M_j\right).$$
Here, the likelihood function $p(\DD |\theta, M_j)$ is replaced by $e^{-\lambda R_n(\DD|\theta, M_j)}$, where $R_n$ is a empirical risk function of the data under the model $M_j$ and the parameter $\theta$. The scaling parameter $\lambda>0$ can depend on the sample size $n$, which is analogous to the inverse temperature in statistical physics. Typically $\lambda\propto n$, as in the usual Bayesian posterior where $-\lambda R_n(\theta)$ is the log likelihood function. However in general we allow $\lambda$ to increase with $n$ at any rate.

This quasi-posterior framework is very useful since it does not need to make as much assumptions on the data generation mechanism as is needed to have a true likelihood function. Although the quasi-posterior originates from other fields such as machine learning and econometrics, research on quasi-posterior from statisticians has been increasing in recent years. It has been applied to problems such as partial likelihood in Cox regression, model-free clustering (\citealt*{bhw16}), and clinically important difference (\citealt{sm17}). The latter involves an interesting case of quasi-posterior with general polynomial convergence rates. The current paper will give two more applications of quasi-posteriors, one incorporating model averaging to cube-root asymptotics, another allowing partial identification.

\subsection{Why our study may be of interest to Bayesian statistics}

Since the Bayesian oracle property is a desirable property for BMA, one naturally hopes that it holds and would like to prove it for some well-established or new methods (see e.g., \citealt{lj16} for Bayesian generalized method of moments, \citealt{ir11} for spike and lab linear regression). Our current paper shows that it is widely valid in the regular cases  for general quasi-posteriors, as long as the model selection consistency holds. This will be useful for Bayesian researchers who invent a new method and would like to go one step further to provide a theoretical justification, in addition to empirical results.

Under the quasi-posterior framework, the more interesting cases are those non-regular models, in which the extremum estimators related to $R_n$ may have nonstandard convergence rates, or the parameters are only partially identified. In such situations, we will show that the Bayesian oracle property does not always hold, and its most straightforward definition may not be always useful. Precaution is needed on how to define a useful oracle property, on how to choose the complexity penalty in the prior, on how to choose the inverse temperature of the quasi-posterior, and on how to choose between BMA and BMS.
From the two examples we study, we find that the answers to the aforementioned questions are highly nontrivial, which could be of interest to Bayesian statisticians.

\subsubsection*{Example 1. Cubic-root asymptotics}

Let $Y=I(Z>0)$ be an observed binary response variable with a latent variable $Z$ related to the utility of the binary choice between $Y=0$ and $Y=1$, where $I(\cdot)$ denotes the indicator function. $Z$ can be modeled as a linear combination of an observed vector of predictors $X$. Given an i.i.d. sample $\DD=\left\{(Y_i,X_i):i=1,\ldots,n\right\}$, one can minimize the empirical risk $R_n(\theta)= -n^{-1}\sum_{i=1}^n Y_iI(X_i^\top \theta\geq 0)$. \citet{m75} discovers that the minimization of $R_n(\theta)$ leads to consistent estimation of $\theta$, when the median of $Z$ is proportional to $X^\top \theta$, without any other distributional assumption on $Z$ such as being normal or logistic. This motivates research on quasi-posteriors (e.g. \citealt*{jpw15}) using $e^{-\lambda R_n(\theta)}$ to play the role of the likelihood function, whose posterior means consistently estimate $\theta$ in a robust way, without additional distributional assumptions on the data. The exponent function $-\lambda R_n(\theta)$ in this example is discontinuous and its minimizers can converge at a rate of $n^{-1/3}$. This is just one example of many similar cases where cubic-root asymptotics appear.

Our study on BMA allows models with various subsets of $X$ components and proves the oracle property, where the asymptotic behavior of the quasi-posterior from BMA is the same as if the true subset of $X$ components were known. In particular, our study in Section 4 shows several nontrivial results in the presence of cubic-root asymptotics:

\noindent 1. {\it On choice of inverse temperature $\lambda$:} The standard choice of $\lambda$ in the likelihoods of regular models is not very useful since it causes the limiting distribution of the posterior mean to be a nonstandard distribution. The BMA has a more useful oracle property when $\lambda$ growers slower than $n^{2/3}$ and faster than $n^{2/5}$.

\noindent 2. {\it  On which oracle property is useful or not useful for the quasi-posterior:}  The oracle property on the quasi-posterior distribution itself is not so useful as a more carefully defined oracle property of the quasi-posterior mean. This is due to the well known result that asymptotically the quasi-posterior distribution may have the correct centering location but the wrong spread. See, e.g., \citet{ch03} who show that the quasi-posterior distributions can give consistent parameter estimates but with wrong standard errors. Therefore, for the purpose of statistical inference, it is more meaningful to consider the mean of the quasi-posterior, rather than the whole quasi-posterior distribution.

\subsubsection*{Example 2. Partial identification}

Consider the example of interval censored data, where an unobservable random variable $Y$ lies in the interval $[L,U]$, and both $L$ and $U$ are observable random variables. The goal is to estimate $\theta=\E (Y)$. Given an i.i.d. sample $\DD=\left\{(L_i,U_i): i=1,\ldots,n\right\}$, one can use the risk function $R_n(\theta)=[(\bar U-\theta)_+]^2+[(\theta-\bar L)_+]^2$ (\citealt*{cht07}), where $(a)_+\equiv \max\{a,0\}$, and $\bar L$ and $\bar U$ are sample averages of $L$ and $U$. The minimizer of $R_n(\theta)$ can be the entire non-singleton set $[\bar L,\bar U]$. A quasi-posterior approach based on this $R_n$ is studied in \citet{w13}. If there exist different prior beliefs in the location of $\theta$, then one can further perform BMA over these different models. A different approach is provided in the example of Section 5, where we use the framework in \citet{ms12} with a reduced-form parameter and a structural parameter. How to properly define BMA and BMS in such partially identified models is very subtle. Through   our effort in finding suitable definitions of Bayesian oracle properties and finding the conditions for them to hold, we obtain several nontrivial results in Section 5 and in a supplementary material, which we believe are of interest to Bayesian statisticians:

\noindent 1. {\it On the formulation of Bayesian oracle properties:} The ``true" model needs to be carefully defined. Partial identification can lead to multiple models that achieve the same minimal risk and are qualified to be the ``true model'' simultaneously. In our simple example above, any model that assigns a uniform prior for $\theta$ in a closed interval can minimize $R_n$ to be zero, as long as this closed interval has non-empty intersection with $[\bar L, \bar U]$. Therefore, it makes more sense to group all such minimum-risk models to form a combined true model in the definition of Bayesian oracle properties, instead of defining the true model as the minimum-risk model with the lowest model complexity.

\noindent 2. {\it On prior choice of complexity penalty:} In the partial identification problem, it is not wise to artificially penalize the model complexity in the prior, in order to favor the simplest minimum-risk model and make it the unique large sample limit in the posterior. In the simple interval censoring example above, suppose that $\E(L)=-0.1$, $\E(U)=0.3$, and the true parameter is $\theta^*=\E(Y)=0.1$. Suppose that one model is given by $\theta\in\{0\}$, i.e. it proposes a singleton prior at $\theta=0$, while the other models do not propose singletons. Then this singleton model achieves the minimum risk zero for $R_n$ asymptotically since $0\in [\E(L),\E(U)]$, but it gives the wrong parameter value since $\E(Y)\neq 0$. Therefore, any penalization through the model priors to favor this simplest but wrong model could lead to misleading inference from the quasi-posterior distribution.

\noindent 3. {\it On BMA versus BMS:} In the presence of partial identification, the oracle property does not hold for the BMS in general. The BMS picks only one of the possibly many minimum-risk models, which may miss the true parameter, as already explained in our first point before. Hence, BMS is not so reliable as BMA, whose limiting quasi-posterior distribution usually includes all those minimum-risk models compatible with the observed data.

In addition to these qualitative guidances on practice, our study also has a number of virtues in theoretical contribution, which are summarized in a technical report \citet{jl15}.

\subsection{Related works}

Bayesian oracle property under model averaging has been considered in the linear model setup by \citet{ir11} and \citet*{csv15}. In contrast, our paper is more general in the sense that it does not assume linear models. \citet{hp12} addressed post selection prediction with possibly nonnested models. \citet{lj16} considered Bayesian generalized method of moments with increasing dimensionality. However, both works assume a regular asymptotic behavior with identifiability and $\sqrt{n}$ asymptotics. The current paper, on the other hand, allows partial identification and cubic-root asymptotics, which entails nonstandard limiting posterior distributions.

We also note that the relationships studied by \citet{hp12} are somewhat different from ours: they relate the point prediction from BMA to the frequentist post-selection predictor, while we study the total variation distance between the entire distributions of the BMA posterior and the oracle posterior given the true model. In this sense, their work and our work are complementary to each other from different perspectives.

\subsection{Organization of the paper}

The rest of the paper is organized as follows. In Section 2 we introduce three types of Bayesian oracle properties for Bayesian model averaging, MAP model selection, and the posterior mean. Section 3 outlines how one can achieve these Bayesian oracle properties in a general quasi-Bayesian framework. These general approaches are then applied to the examples of cubic root asymptotics in Section 4 and partially identified models in Section 5. Section 6 summarizes the paper with some discussions. Section 7 contains the proofs of the propositions. All other technical details and proofs are included in a supplementary material.

We introduce some useful notation. For two $n$-dependent sequences $\{a_n\}$ and $\{b_n\}$, $a_n\prec b_n$ and $b_n\succ a_n$ denote the relation $\lim_{n\to\infty} a_n/b_n=0$. $a_n\preceq b_n$ and $b_n\succeq a_n$ denote that $a_n/b_n$ is bounded by constant. $a_n\asymp b_n$ is equivalent to $a_n\preceq b_n$ and $b_n\preceq a_n$. We use $I(\cdot)$ to denote the indicator function. We use $o_p(1)$ and $O_p(1)$ to denote the orders under the probability measure of $\mathbf{D}$ as the sample size $n$ increases to infinity.

\section{Bayesian Oracle Properties}
\subsection{Bayesian model averaging}

The first property we define here is the global model selection consistency.

\begin{property}\label{O1}
$\pi(\theta|\DD)$ satisfies the global model selection consistency, if $1-\pi(M^*|\DD)=o_p(1)$.
\end{property}

The global model selection consistency says that the true model $M^*$ has posterior probability converging to 1 as the sample size increases to infinity. The consistency holds for the regular parametric model under the Bayesian framework, based on the standard BIC theory (\citealt{s78}). It also holds for general high dimensional regression models under certain priors that induce sparsity (\citealt{jr12}, \citealt{lsy13}, etc.).

For any (data-dependent) measurable event $A$, we are interested in the difference between two probabilities
$$\left|\Pi(A|\DD)-\Pi(A|M^*,\DD)\right|,$$
where $\Pi(A|M^*,\DD)$ is the probability of $A$ under the ``oracle" posterior distribution, pretending that the true model $M^*$ is known, whereas $\Pi(A|\DD)=\sum_{j\geq 1} \pi(M_j|\DD)\cdot \Pi(A|M_j,\DD)$ is the mixed posterior distribution via model averaging, allowing possibilities of all models which are weighted by the model posterior probabilities $\pi(M_j|\DD)$ for $j=1,2,\ldots$.

\begin{property}\label{O2}
$\pi(\theta|\DD)$ satisfies the oracle property for Bayesian model averaging, if\\
$\sup_{A\in \cal {F}}\left|\Pi(A|\DD)-\Pi(A|M^*,\DD)\right|=o_p(1)$ where $\cal {F}$ is the set of all measurable events.
\end{property}

This defines an oracle property for Bayesian model averaging, which basically says that any posterior inference based on model averaging is asymptotically equivalent to the oracle posterior inference based on only the true model. It turns out that one can establish the following fundamental inequality.
\begin{prop}\label{tvmisprob}
$$ \sup_{A\in\cal {F}}\left|\Pi(A|\DD)-\Pi(A|M^*,\DD)\right| \leq  1-\pi(M^*|\DD),$$
where $\cal {F}$ is the set of all measurable events.
\end{prop}

This proposition reveals a deep relation between three quantities: the model averaging posterior $\pi(\theta|\DD)$, the oracle posterior $\pi(\theta|M^*,\DD)$, and the posterior probability of the true model $\pi(M^*|\DD)$. The total variation distance between the model averaging posterior and the oracle posterior is bounded above by the posterior probability of missing the true model. A direct consequence of Proposition \ref{tvmisprob} is the relation between the global model selection consistency (Property \ref{O1}) and the oracle property for Bayesian model averaging (Property \ref{O2}).

\begin{thm}
The global model selection consistency (Property \ref{O1}) implies the oracle property for Bayesian model averaging (Property \ref{O2}).
\end{thm}

Therefore, as the sample size increases to infinity, if the true model has posterior probability converging to 1, then the limiting behavior of the posterior distribution under model averaging is the same in total variation norm as the oracle posterior pretending to have known the true model. This kind of oracle property is similar in essence to the frequentist oracle property of \citet{fl01} but is more general.

To fully appreciate the generality of Theorem 1, we emphasize that in the current general context, we do not require the oracle posterior $\pi(\theta|M^*,\DD)$ to satisfy the parametric Bernstein von Mises theorem (BvM), i.e. converging to a normal limiting distribution asymptotically at the rate of $n^{-1/2}$. The most attractive aspect of Fan and Li's oracle property is that the inferential results ``work as well as if the correct submodel were known" (see the abstract of \citealt{fl01}). This aspect has already been fully captured by Property \ref{O2} and there is no need to impose any additional restrictions on the oracle posterior $\pi(\theta|M^*,\DD)$. Our relaxation makes it possible to include many nonstandard models where a parametric BvM type result does not hold, such as the (quasi-)posteriors with discontinuous (quasi-)likelihoods which is characterized by the cubic root asymptotics (see, e.g., \citealt*{jpw15}), and the partially-identifying posterior distributions with the $O(1)$ rate asymptotics (see, e.g., \citealt{ms12}).

\subsection{MAP (maximum a posteriori)  model selection}\label{sec.select}
As an alternative to Bayesian model averaging, one could select only one MAP model that has the maximum posterior probability. We would like to establish  similar results to Theorem 1 for MAP model selection. Suppose $\widehat M$ is any MAP model choice, so that $\pi(\widehat M|\DD)=\max_{j\geq 1} \pi(M_j|\DD)$. We are interested in the total variation distance between the posterior $\pi(\theta|\widehat M, \DD)$ based on the MAP model, and the oracle posterior  $\pi(\theta|M^*,\DD)$ based on the true model $M^*$. We hope that inference based on the MAP model choice $\widehat M$ is almost as good as if based on the true model $M^*$.

\begin{property}\label{O3}
$\pi(\theta|\DD)$ satisfies the oracle property for MAP model selection, if \\
$\sup_{A\in{\cal F}}\left|\Pi(A|\widehat M,\DD)-\Pi(A|M^*,\DD)\right|=o_p(1)$ where $\cal {F}$ is the set of all measurable events.
\end{property}

Based on this definition, we have the following proposition.

\begin{prop}\label{tv3}
The maximal total variation distances among any of the three posteriors
$\Pi(\cdot|\DD)$, $\Pi(\cdot|\widehat M,\DD)$, and $\Pi(\cdot|M^*,\DD)$,
are at most twice the posterior probability of missing the true model $2\left[1-\pi(M^*|\DD)\right]$.
\end{prop}

A direct consequence of this proposition is
\begin{thm}
The global model selection consistency (Property \ref{O1}) implies the oracle property for MAP model selection (Property \ref{O3}).
\end{thm}

\subsection{Mean oracle property}

In some situations the (quasi-)posterior $\pi(\theta|\DD)$ itself is either not of main interest or does not have any valid interpretation, but the posterior mean $\E(\theta|\DD)=\int_{\Theta}\theta d\pi(\theta|\DD)$ for some parameter $\theta$ is still of interest, which may have a well understood limiting distribution that can be used for inference on $\theta$. This can happen for quasi-posteriors when its credible region does not have asymptotically correct coverage probability. One example is the Bayesian quantile regression with a quasi-likelihood constructed from the check function. The generalized information criterion is violated and the quasi-posterior has no valid interpretation (\citealt{ch03}), but the posterior mean can be used as a convenient frequentist estimator for the quantile regression coefficients. Another example is the Laplace version of the least median of squares estimator (\citealt*{jpw11}). In this case, it is desirable to have a version of {\em Bayesian oracle property for the posterior mean: If we make inference based on the overall posterior mean, it is as if we were making inference based on the posterior mean conditional on the true model only}.

To achieve such oracle inference for the mean for a posterior distribution $\pi(\cdot)$, it is usually not sufficient to only have the relation $\|\E(\theta|\DD)-\E(\theta|M^*,\DD)\|=o_p(1)$, because $\E(\theta|\DD)$ and $\E(\theta|M^*,\DD)$ may both converge to a true parameter $\theta^*$ but with different convergence rates. A more proper version of mean oracle property is defined as follows.

\begin{property}\label{O4}
$\pi(\theta|\DD)$ satisfies the mean oracle property, if\\
$\|\E(\theta|\DD)-\E(\theta|M^*,\DD)\|= o_p(1) \cdot \|\E(\theta|M^*,\DD)-\theta^*\|.$
\end{property}
In other words, we require that the difference between posterior means from Bayesian model averaging and the oracle is of higher order compared to the posterior bias under the oracle posterior. This will guarantee that $\E(\theta|\DD)-\theta^*$ and $\E(\theta|M^*,\DD)-\theta^*$ are approximately the same, and not merely both converging to zero.

A useful relation which can be applied to achieve the mean oracle property is
\begin{align} \label{detailrate}
& \E(\theta|\DD)-\E(\theta|M^*,\DD) = \sum_{j\geq 1, M_j\neq M^*}\pi(M_j|\DD)\left[\E(\theta|M_j,\DD)-\E(\theta|M^*,\DD)\right].
\end{align}
The mean oracle property holds if there is a fixed number of model candidates and for every model $M_j\neq M^*$,
$\pi(M_j|\DD)\left\|\E(\theta|M_j,\DD)-\E(\theta|M^*,\DD)\right\| = o_p(1) \left\|\E(\theta|M^*,\DD)-\theta^*\right\|$. Each product in the sum of (\ref{detailrate}) can be made small enough for different reasons. For example, consider the standard variable selection problem in linear models. For those models that miss nonzero parameters, $\pi(M_j|\DD)$ is typically exponentially small.
For the models that do not miss nonzero parameters but include redundant parameters, $\E(\theta|M_j,\DD)-\theta^*$ is typically of the same order as $\E(\theta|M^*,\DD)-\theta^*$, and therefore $\E(\theta|M_j,\DD)-\E(\theta|M^*,\DD)$ is also of the same order as $\E(\theta|M^*,\DD)-\theta^*$; then it is sufficient to have $\pi(M_j|\DD)=o_p(1)$.
The method described here will be applied to a nonstandard example with cubic-root asymptotics in Section \ref{cubicrootasymp}.

\subsection{Applications}
There has been extensive work in Bayesian model selection consistency, especially the global model selection consistency (Property \ref{O1}). All these results can be readily extended to the oracle property for Bayesian model averaging (Property \ref{O2}) and for MAP model selection (Property \ref{O3}). Whenever there are already known results on the limiting distribution of the oracle posterior $\pi(\theta|M^*,\DD)$ under the true model $M^*$, the limiting distribution automatically applies to $\pi(\theta|\DD)$ from model averaging by Theorem 1 and to $\pi(\theta|\widehat M,\DD)$ from model selection by Theorem 2.

The most well known example is the regular finite dimensional models, where BvM type results hold and the posterior distribution of finite dimensional parameters converges in total variation norm to the normal limit at the parametric rate of $n^{-1/2}$. See for example, Section 10.2 of \citet{van98} for finite dimensional parametric models, and \citet{s02} for nonparametric and semiparametric models. Consequently, in combination with the classic BIC theory from \citet{s78}, one can derive the global model selection consistency (Property \ref{O1}) for such finite dimensional cases (see for example \citealt{w00} Equation 42), and our theorems suggest that the posterior inference based on model averaging or model selection  is also equivalent to the inference under the limiting normal distribution given the (unknown) true model. When the model is regular and high dimensional, exactly the same equivalence holds as long as a BvM type result can be established for the low dimensional true model $M^*$, with properly chosen sparsity inducing priors, such as the priors used in \citet{jr12} and \citet{lsy13}.

In this paper, we are interested in applications of the Bayesian oracle property under a more general Bayesian framework  than the regular parametric models. We extend the likelihood-based posterior to the general quasi-posterior, in which the likelihood function is replaced by a quasi-likelihood based on a risk function. We propose two ways to achieve the Bayesian oracle properties in Section \ref{bic01} and \ref{bic02} respectively, with two applications: The first application is to the cubic root asymptotics where the convergence rate is not the standard parametric rate $n^{-1/2}$. The second application is to partially identified models where the posterior distribution has a nonstandard limit and a BvM type result does not hold.

\section{Quasi-posterior with General Risk}\label{quasipost}
We will work under the general framework of a (quasi-)posterior where we can derive general bounds on the mis-selection probability $1-\pi(M^*|\DD)$. As discussed in Section 1.2, we consider the quasi-posterior distribution
\begin{align}\label{quasipost0}
&\pi(\theta|\DD) = \frac{ e^{-\lambda R_n(\theta)} d\pi(\theta) }{\int_{\Theta} e^{-\lambda R_n(\theta) } d\pi(\theta)},
\end{align}
where $\pi(\theta)$ is the prior density and $R_n$ is an empirical risk function dependent on both the parameter $\theta$ and the data $\DD$.  Related to $R_n(\theta)$ is a theoretical risk function $R(\theta)$, which is typically the large sample limit of $R_n(\theta)$. The scaling parameter $\lambda>0$ can depend on $n$ and increase with $n$ at any rate, which is analogous to the inverse temperature in statistical physics.

We describe what a true model and a true parameter mean. This is not always clear in the context of quasi-posteriors. Since our quasi-posterior is related to an empirical risk $R_n(\theta)$, which usually has a theoretical risk $R(\theta)$ as its large sample limit, we will treat the minimizer of $R(\theta)$ over the entire parameter space $\Theta$ as our true parameter $\theta^*$. We will define a minimum-risk model to be {\it a model whose prior support includes $\theta^*$}. Situations can be complicated in that there may be multiple minimum-risk models. Conventional wisdom suggests defining the true model $M^*$ as the simplest minimum-risk model that has the lowest dimension of the prior support. If needed, we can also group multiple minimum-risk models together as a composite true model with a mixture prior. A later Section \ref{partial_identify} uses this approach to handle partial identification, where the minimizer of $R(\theta)$ is not a singleton and some variation is needed in defining the true model.

In the following, we consider two methods of bounding $1-\pi(M^*|\DD)$, the quasi-posterior probability of mis-selecting the true model. Our results from previous sections have shown that bounding this mis-selection probability can lead to various oracle properties. We will make an assumption of finitely many models for simplicity.


\subsection{Bounding the mis-selection probability: Extending the BIC approximation for quasi-posterior} \label{bic01}

In the classical BIC approach (\citealt{s78}), a complexity penalty arises indirectly from approximating an integral in the posterior calculation. Suppose that the parameter space $\Theta_j$ is finite dimensional for any $j\geq 1$ and the dimension $d_j=\dim(\Theta_j)$ is bounded. Let $\Theta^*$ be the parameter space of $M^*$ and $d^*=\dim(\Theta^*)$. The prior probabilities $\pi(M_j)$ are all assumed to be of order 1 and will not affect the asymptotic behavior. Suppose that the risk functions $R(\theta)$ and $R_n(\theta)$ only depend on the value of $\theta$ and do not depend on the model index $M_j$. For convenience we assume that $\theta^*=\arg\min_{\theta\in \Theta} R(\theta)$ is the unique minimizer of $R(\theta)$. We can extend the BIC approximation to general quasi-posteriors and bound the posterior mis-selection probability.

\begin{prop} \label{bicproperty}
Consider the following assumptions: \\
\noindent (i) The total number of models is bounded above by a constant integer, and all models have a positive prior probability; \\
\noindent (ii) For any minimum-risk model $M_j$ that satisfies $\inf_{\Theta_j} R(\theta) = R(\theta^*)$ (which implies $\theta^*\in \Theta_j$), the integral in the posterior model probability satisfies a BIC type approximation
\begin{align}\label{bic11}
-\ln \int_{\Theta_j} e^{-\lambda R_n(\theta)} d\pi(\theta|M_j)
=  \lambda R_n(\theta^*) + \frac{d_j\ln\lambda}{2} + O_p(1);
\end{align}
\noindent (iii) For any minimum-risk model $M_j\neq M^*$, $d_j\geq d^*+1$; \\
\noindent (iv) For any non-minimum-risk model $M_j$ with $\inf_{\Theta_j} R(\theta) - R(\theta^*) \equiv \gamma_j>0$, we have $\gamma_j\succeq 1$ and $S_n(\theta)=o_p(1/\lambda)$ uniformly over $\theta\in \Theta_j$, where $S_n(\theta)=[R_n(\theta)-R(\theta)]-[R_n(\theta^*)-R(\theta^*)]$;\\
\noindent (v) $\lambda\to \infty$ as $n\to\infty$; \\
\noindent (vi) $\Theta$ is compact. The scaling parameter $\lambda$ grows polynomially in $n$. For any minimum-risk model $M_j$,
$$\left\|\E(\theta|M_j,\DD)-\theta^*\right\|=O_p(1)\cdot \left\|\E(\theta|M^*,\DD)-\theta^*\right\|,$$
and $\left\|\E(\theta|M^*,\DD)-\theta^*\right\|\succeq \epsilon_n$, where $\epsilon_n=o(1)$ is polynomial in $n$.\\
\vspace{-5mm}

Then Bayesian oracle properties \ref{O1}, \ref{O2}, and \ref{O3} hold under the assumptions (i)-(v), and the Bayesian mean oracle property \ref{O4} holds under the assumptions (i)-(vi).
\end{prop}

Although the approach outlined in this subsection is still mathematically a BIC approximation, it is somewhat more general, in that it accommodates non-likelihood based quasi-posterior and an arbitrary scaling $\lambda$ that may increase at a different rate than $n$. It turns out that this extension of BIC can be applied to the example with nonstandard cubic-root asymptotics in Section \ref{cubicrootasymp}.


%

\subsection{Bounding the mis-selection probability: Assumption-free upperbound for quasi-posterior} \label{bic02}

In the later example with partial identification (Section \ref{partial_identify}), the BIC approximation (which uses a local approximation of the theoretical risk $R$ near its minimum) will no longer work. We will apply the following assumption-free upper bound on the mis-selection probability $1-\pi(M^*|\DD)$, which does not require $\arg\min_{\theta\in \Theta} R(\theta)$
to be a singleton, and can therefore be applied to situations with partial identification.

\begin{prop}\label{gibbs1}
(Model selection with quasi-posterior)
The mis-selection probability is upper bounded by
$$\ln [1-\pi(M^*|\DD)] \leq -0.5\lambda (\gamma -r-2|u|)$$
where
\begin{align*}
\gamma& =\inf_{\theta\in \Theta, M\neq M^*} R(\theta) - \inf_{\theta\in \Theta} R(\theta),\\
r &=-\lambda^{-1}\ln \int_{\Theta} e^{-\lambda[R(\theta)- \inf_{\theta\in \Theta} R(\theta)]}\pi(\theta)d\theta, \\
u &=-(2\lambda)^{-1}\ln \int_{\Theta} e^{-2\lambda\left[(R_n(\theta)-R(\theta))-\int_{\theta\in \Theta} (R_n(\theta)-R(\theta))\pi_{\infty}(\theta)d\theta\right]} \pi_{\infty}(\theta) d\theta,
\end{align*}
and
\begin{align*}
\pi_{\infty}(\theta) =\frac{e^{-\lambda R(\theta)}\pi(\theta)d\theta} {\int_{\theta\in \Theta} e^{-\lambda R(\theta)}\pi(\theta)d\theta }
\end{align*}
is the limiting version of the quasi-posterior $\pi(\theta|\DD)$, in which the theoretical risk $R$ is used in place of the empirical risk $R_n$.
\end{prop}

This assumption-free bound uses three quantities: $\gamma$ ($gap$), which differentiates the best possible
risks achievable by model $M^*$ and by other models; and $r$ ($excess$), which is a nonstochastic term related to the excess risk $R(\theta)-\inf_{\theta\in \Theta} R(\theta)$ which we will bound later; $|u|$ ($noise$), which is a stochastic noise term determined by the difference $R_n(\theta)-R(\theta)$.
This assumption-free bound is only useful when $\gamma > r+2|u| > 0$. We show in the following how it is possible to make $r+2|u| = o_p(\gamma)$, such that $1-\pi(M^*|\DD)$ can be exponentially small in $\lambda \gamma$ and decreases very quickly with sample size $n$.

The noise term $u$ measures the difference $R_n(\theta)-R(\theta)$ on the support of the limiting posterior. We can use the simplest uniform bound
$$ |u|\leq 2\sup_{\theta\in \Theta}\left|R_n(\theta)-R(\theta)\right|.$$
By using uniform large deviation, this will typically lead to $u=O_p(\ln n/\sqrt{n})$. The nonstochastic term $r$ can be bounded by $r=O(\ln\lambda /\lambda)$ if $R(\theta)$ allows a Laplace approximation.

In general, without assuming a Laplace approximation for $R(\theta)$, the rate $r=O(\ln\lambda /\lambda)$ can be derived by the inequality
\begin{align}\label{ljt}
r &= -\lambda^{-1}\ln \int_{\theta\in \Theta} e^{-\lambda [R(\theta)-\inf_{\theta\in \Theta} R(\theta)]} \pi(\theta)d\theta \nonumber \\
&\leq \inf_{a>0} \left[a - \frac{1}{\lambda} \ln \Pi\left(\left\{\theta:R(\theta)-\inf_{\theta\in \Theta} R(\theta)<a\right\}\right)\right],
\end{align}
and choosing $a=\ln\lambda /\lambda$. Detailed argument is similar to the remarks after Proposition 1 in \citet*{ljt14}.

Therefore, if $\gamma\succ r+2|u|$ and $\gamma\succ \ln n/\lambda$, then $1-\pi(M^*|\DD)=\pi(M\neq M^*|\DD)\prec e^{-\ln n}=1/n\to 0$ as $n\to\infty$, and we achieve the global model selection consistency (Property \ref{O1}). Therefore the oracle properties \ref{O2} and \ref{O3} also hold true.
The above bound for $1-\pi(M^*|\DD)$ may also be used to prove the mean oracle property \ref{O4} with the help of \eqref{detailrate}.

\section{Cubic Root Asymptotics}\label{cubicrootasymp}

Suppose that we observe i.i.d. data $\mathbf{D}=\left\{D_1,\ldots,D_n\right\}$, and the parameter of interest is $\theta\in \Theta \subseteq \mathbb{R}^p$, whose true value $\theta^*$ is the unique solution to the optimization problem $\min_{\theta \in \Theta }Eg(D_1,\theta)$ for some known criterion function $g$ and the expectation is taken with respect to the true underlying distribution of $D_1$. Let $R(\theta)=\E g(D_1,\theta)$ be the theoretical risk and $R_n(\theta)= n^{-1}\sum_{i=1}^n g_i(\theta)$ be the empirical risk where $g_i(\theta)$ is a shorthand for $g(D_i,\theta)$. Instead of the parametric rate $n^{-1/2}$, the frequentist extremum estimator which minimizes $R_n(\theta)$ may have a slower $n^{-1/3}$ convergence rate when $g$ is discontinuous in $\theta$. For example, if one predicts a binary variable $Y_i$ with a vector of continuous predictors $(X_{0,i},X_i)^\top\in \mathbb{R}^{p+1}$, the maximum score estimator (\citealt{m75}) minimizes $R_n(\theta)$ with $g_i(\theta)=-Y_iI(X_i^\top \theta - X_{0,i} \geq 0)$, which asymptotically can have a $n^{-1/3}$ convergence rate. Here we assume that the variable $X_{0,i}$ is always selected and its coefficient is $-1$ to ensure the identification of $\theta^*$.  Other applications of the cubic root asymptotics include shorth estimation, least median of squares estimator, isotonic regression, quantile regression with interval censoring, etc. See \citet{kp90} and \citet*{jpw15} for more examples.

We consider the quasi-posterior defined in \eqref{quasipost0} using the empirical risk function $R_n(\theta)=n^{-1}\sum_{i=1}^n g_i(\theta)$. For the ease of presentation, we only consider the ``theta class" in \citet*{jpw15}. The Laplace type estimator of $\theta$ discussed in \citet*{jpw15} is the posterior mean of \eqref{quasipost0}. The standard model/variable selection in this cubic root problem assumes that the true parameter $\theta^*$ could possibly lie in a lower dimensional space $\Theta \cap \mathbb{R}^{p^*}$ with $1\leq p^*\leq p$. For example, for the maximum score estimator, our goal is to select only the relevant predictors in $X$ and we set the $\theta$ coefficients of all irrelevant components of $X$ to be zero. Then a model $M_j$ in this context is defined as a coordinate subspace of $\Theta\cap \mathbb{R}^{p}$. The maximum number of possible models in $\Theta\cap \mathbb{R}^p$ is $2^p-1$. The true model $M^*$ is defined to be the lowest dimensional coordinate subspace that contains the true parameter $\theta^*$ such that all components of $\theta^*$ in $M^*$ are nonzero. We assume that the prior density has the decomposition $\pi(\theta|M_j)\pi(M_j)$, where $\pi(\theta|M_j)$ is a continuous density on $\Theta_j \equiv \Theta \cap \mathbb{R}^{d_j}$, $d_j$ is the dimension of $M_j$, and $\sum_{j=1}^{2^p-1}\pi(M_j)=1$ give the discrete probabilities for all models.

We make the following assumptions on the model and the prior.

\begin{itemize}
\item[(C1)] $\Theta$ is compact. $\theta^*$ is an interior point of $\Theta\subseteq \mathbb{R}^p$ with $p$ being a constant dimension. $\min_{j\in M^*}|\theta^*_j|\geq c_{\theta}$, where $\theta^*_j$ for $j\in M^*$ denotes the $j$th nonzero component of $\theta^*$ and $c_{\theta}>0$ is a constant.
\item[(C2)] For all $\theta \neq \theta^*$, $R(\theta)>R(\theta^*)$.
\item[(C3)] $R(\theta)= \E g(D_1,\theta)$ is three times continuously differentiable in $\Theta$. Let $V=\partial_{\theta\theta^\top} R(\theta^*)$ be the second derivative matrix of $R(\theta)$ evaluated at $\theta=\theta^*$. Then $V$ is positive definite with eigenvalues bounded from below and above by positive constants.
\item[(C4)] For any $t,s\in \mathbb{R}^p$, the function $H(t,s)=\lim_{a\to +\infty} a\E \left[g_1(\theta^*+t/a)g_1(\theta^*+s/a)\right]$ exists and is always positive.
\item[(C5)] $\pi(\theta|M_j)$ is continuously differentiable for all $\theta\in \Theta_j$ and all models $M_j$. $\pi(\theta|M_j)$ and $\partial_{\theta} \pi(\theta|M_j)$ are uniformly bounded from above by constant for all $\theta\in \Theta_j$ and all models $M_j$. For all models $M_j$ that satisfy $M_j\supseteq M^*$,  $\pi(\theta^*|M_j)$ is uniformly bounded from below by a positive constant. $\pi(M_j)$ is bounded from above and below by positive constants for all models $M_j$.
\end{itemize}

Similar to \citet*{jpw15}, we make the following assumptions on the envelope function of $g(D_1,\theta)$. These assumptions depend on the inverse temperature parameter $\lambda$ in the quasi-posterior \eqref{quasipost0}. Let $g^\circ(D_1,t)=\lambda^{1/4}[g(D_1,\theta^*+t/\sqrt{\lambda})-g(D_1,\theta^*)]$ $/(\|t\|+1)$. Let $\mathcal{G}_n=\left\{g^\circ(D_1,t):t\in \mathbb{R}^p\right\}$.

\begin{itemize}
\item[(C6)] There exists an envelope function $G(\cdot)$ such that $\sup_{t\in \mathbb{R}^p} \left|g^\circ (D_1,t)\right|\leq G(D_1)$ almost surely under the distribution of $D_1$. Furthermore,  $\E\left[G^2(D_1)\right]<\infty$ and $\lim_{n\to\infty}\E [G^2(D_1)\cdot
    I(G(D_1)>c\sqrt{n})]=0$ for any $c>0$.
\item[(C7)] For any $0<\epsilon_n=o(1)$, $\sup_{t,s\in \mathbb{R}^p, \|t-s\|\leq\epsilon_n} \E \left[g^\circ(D_1,t)-g^\circ(D_1,s)\right]^2=o(1).$
\item[(C8)] Let $\mathcal{N}(\epsilon,\mathcal{G}_n,L_2(P))$ be the $L_2$-covering number for $\mathcal{G}_n$ with respect to the probability measure $P$. Then for every sequence  $0<\epsilon_n=o(1)$,
$$\sup_{P^*}\int_0^{\epsilon_n} \sqrt{\log \left[\mathcal{N}\left(\epsilon\|G(D_1)\|_{P^*},\mathcal{G}_n,L_2(P^*)\right)\right]}d\epsilon =o(1),$$
where $\sup_{P^*}$ is the supremum taken over all finitely discrete probability measures $P^*$ with $\|G(D_1)\|_{P^*}=\sqrt{\E_{P^*}[G^2(D_1)]}>0$.
\end{itemize}

(C1) assumes the standard beta-min condition on $\theta^*$ to distinguish its nonzero and zero components. We use the constant lower bound $c_{\theta}$ for technical convenience, as it could be replaced by a rate slowly decreasing to zero that depends on the growth rate of $\lambda$. (C2)-(C4) and (C6)-(C8) are similar to the conditions used in \citet*{jpw15}, which leads to the cubic root behavior of the frequentist extremum estimator that minimizes $R_n(\theta)$. (C5) contains mild conditions on the model selection prior. The essential requirement is that every plausible model should have positive prior probabilities, and the prior mass around the true parameter $\theta^*$ should not be too small.

\begin{thm}\label{cubicroot}
Suppose (C1)-(C8) hold with $\lambda$ satisfying $n^{2/5}\prec \lambda \prec n^{2/3}$. Then the global model selection consistency (Property \ref{O1}), the Bayesian model averaging oracle property (Property \ref{O2}),  the MAP model selection oracle property (Property \ref{O3}),  and the mean oracle property (Property \ref{O4}) all hold for the quasi-posterior $\pi(\theta|\DD)$ in \eqref{quasipost0}.
\end{thm}

In Theorem \ref{cubicroot} we restrict the growth rate of $\lambda$ to be between $n^{2/5}$ and $n^{2/3}$. The main reason is that with such $\lambda$, the limiting distributions of both the quasi-posterior and the posterior mean will be normal with mean zero, even under a model selection setup with our condition (C5) on the prior. The contribution of our mean oracle property basically says that the asymptotics of the posterior mean from \citet*{jpw15}, who did not consider model selection but assumed the true model to be known, still remains valid as if the true model were known when we have a pool of candidate models with an unknown true model.

The conclusion of Theorem \ref{cubicroot} follows from the BIC type approximation in Case (iii) of Theorem 1 in \citet*{jpw15} together with our approach in Section \ref{bic01}. A heuristic argument is as follows. The exponent in the quasi-posterior \eqref{quasipost0} has the decomposition $\lambda R_n(\theta) = \lambda [R(\theta)-R(\theta^*)] + \lambda S_n(\theta)$ with $S_n(\theta)$ defined in Proposition \ref{bicproperty}. Although $R_n(\theta)$ is discontinuous in $\theta$, $R(\theta)$ is continuously differentiable in $\theta$ by (C3). As a result, for any model $M$ that includes the true model $M^*$ as a submodel (including $M^*$ itself), we have a quadratic approximation $\lambda [R(\theta)-R(\theta^*)] \asymp \lambda\|\theta-\theta^*\|^2$. Meanwhile it can be shown that the $S_n(\theta)$ term has a Gaussian process limit and is about the order $O_p(n^{-1/2}\|\theta-\theta^*\|^{1/2})$. Therefore the nonstochastic term of $\lambda [R(\theta)-R(\theta^*)]$ will dominate the stochastic term $\lambda S_n(\theta)$ if $\lambda n^{-1/2}\|\theta-\theta^*\|^{1/2} \prec \lambda\|\theta-\theta^*\|^2 \asymp 1$ in the asymptotics, which leads to $ \lambda \prec n^{2/3}$ and $\|\theta-\theta^*\|\asymp \lambda^{-1/2}$. Hence the BIC approximation in Proposition \ref{bicproperty} works for the minimum-risk models. For any wrong model $M$ that misses at least one component of $M^*$, it follows from the aforementioned relations that $S_n(\theta)\asymp n^{-1/2}\lambda^{-1/2}\prec 1/\lambda$, which implies that the integral in \eqref{bic12} is $O_p(1)$. Hence these models will have exponentially small posterior probabilities in $\lambda$.
The other condition $\lambda \succ n^{2/5}$ in Theorem \ref{cubicroot} is required to eliminate the asymptotic bias of the posterior mean. See the comments after Theorem 1 of \citet*{jpw15}. As a result, the global model selection consistency and the Bayesian oracle properties (Properties \ref{O1}-\ref{O4}) hold true following the argument in Section \ref{bic01}.

The slowly growing $\lambda$ in Theorem \ref{cubicroot} can overcome the discontinuity in the empirical risk $R_n(\theta)$ with a smoothing effect and justifies the BIC type approximations. The posterior convergence rate is $\lambda^{-1/2}$ from the BIC approximation discussed above, which is slower than $n^{-1/3}$ due to the condition on $\lambda$. The posterior mean has a different convergence rate of $n^{-1/2}\lambda^{1/4}$ (see \citealt*{jpw15} Theorem 1 (iii)), which is faster than $n^{-1/3}$.

In this cubic root example, although the limiting distribution of $\pi(\theta|\DD)$ in \eqref{quasipost0} is normal, the quasi-posterior itself typically does not have the usual Bayesian interpretation even in the asymptotic sense of \citet{ch03}. Therefore, the MAP model selection oracle property (Property \ref{O3}) and the model averaging oracle property (Property \ref{O2}) are not meaningful, since the quasi-Bayesian inference based on the true model may still be invalid. However, the mean oracle property (Property \ref{O2}) can be very useful because the posterior mean can converge faster than $n^{-1/3}$ to a limiting normal distribution, under the choice of $\lambda$ in Theorem \ref{cubicroot}. The normal limit allows us to use various tools such as bootstraps or subsampling to construct asymptotically valid confidence intervals for the posterior mean estimator. Hence statistical inference based on the posterior mean estimator can be more advantageous than that based on the frequentist extremum estimator whose limiting distribution is the Chernoff's distribution (\citealt{kp90}). \\

\section{Partial Identification}\label{partial_identify}

In econometric and statistical literatures, there exist two different approaches to handle partial identification. One aims for more informative inference about the partially identified point parameter $\theta$ by incorporating prior information (see, e.g., \citealt{p98}, \citealt{ms12}, \citealt{gus15}). Another aims for more robust inference about the fully identified identification region $\Omega$ (see, e.g., \citealt{w13}, \citealt{kt16}, and \citealt*{cct16}). The current paper follows the first approach.

In this section, we apply Bayesian model averaging to a situation with partial identification as described in \citet{ms12}, who showed that the limiting posterior is nonstandard. The posterior contraction rate for a structural parameter for interest is typically of order 1, instead of the classical order $n^{-1/2}$,  due to partial identification. For example, Equation (4) of \citet{ms12} provides a simple example where the limiting posterior for the structural parameter of interest is uniform over an non-shrinking interval. Despite such nonstandard limiting behavior with partial identification, our machinery in Section \ref{bic02} (based on bounding the mis-selection probability) can be used to study the oracle properties under Bayesian model averaging, which uses a  conservative approach to preserve all submodels that are compatible with the data.

\subsection{A simple example}\label{simpleexample}

This example is similar to the simple example in \citet{ms12}. We add the aspect of model selection or model averaging, and make a small variation that a quasi-likelihood is used instead of a real likelihood.
Suppose we are interested in a structural parameter $\omega=\E Y$, where $Y\in[0,4]$ is the GPA of a college student. However, the GPA is sometimes only known to fall in some interval. For simplicity, assume only its integer part $Z=\lfloor Y\rfloor$ of the GPA is observed. The fractional part $U=Y-Z$ is unobserved. We define $\E Z=\phi$, which is called the reduced-form parameter which is identified by the observed data $Z$. We will call the ``combined'' parameter $\theta=(\omega,\phi)$.  Note that $Z\in [Y-1, Y]$, and therefore $\phi=\E Z\in [\E Y-1, \E Y]=[\omega-1,\omega]$.

In Bayesian approach, the relation between $\phi$ and $\omega$ is described by a conditional prior distribution $\pi(\phi|\omega)$ such as $\text{const}\times I(\{\phi\in [\omega-1,\omega]\cap[0,4]\})$. This conditional prior will be assumed to be the same for all models that we will consider, since we are interested in model selection or model averaging on the structural parameter $\omega$ only. Each candidate models $M_j$, indexed by $j=1,2,\ldots$ and weighted by $\pi(M_j)$, proposes a different prior $\pi(\omega|M_j)$ for the structural parameter $\omega$. So the joint prior for the combined parameter $\theta$ and $M_j$ is $$\pi(\theta,M_j)=\pi(M_j)\pi(\omega|M_j)\pi(\phi|\omega).$$
This way, we  can convert the model selection  problem for the structural parameter $\omega$ to a model selection problem with the combined parameter $\theta$. This is for a technical reason to apply the framework of Section \ref{quasipost} in establishing the oracle properties with Bayesian model averaging, later in Section \ref{oraclepid}.

We will introduce some related concepts first for a very simple example, where $j=1,2$, $\pi(M_j)=1/2$,
$\pi(\omega|M_1)=\delta_{3}(\omega)$ is a point mass supported on $W_1=\{3\}$, proposing mean GPA to be 3, and $\pi(\omega|M_2)=0.25I(\{\omega\in [0,4]\})$ is a prior supported on $W_2=[0,4]$, proposing no restriction on the mean GPA. This can be regarded as a simplified version of the example in the supplementary material, where Figure 1 illustrates prior densities for more than two candidate
models, the first two of them being the same as the current models with $j=1,2$.

The observed data $Z$ is integer valued and nonnormal. However, we can use a normal $quasi$-likelihood based on $\bar Z$, (the observed sample average of $Z$), which is typically asymptotically normal iid data: $\sqrt{n/\hat v}(\bar Z-\phi)\rightarrow N(0,1)$ as $n\rightarrow\infty$, where $n$ is the sample size, and $\hat v$ is a consistent estimate of $v={\rm var}(Z)$. Then the corresponding quasi-posterior has the form $\pi(\theta,M_j)\propto e^{-\lambda R_n(\theta)} \pi(\theta,M_j)$, where $\lambda=n$ and $R_n(\theta)=0.5\hat v^{-1} (\bar Z-\phi)^2$ is an empirical risk derived from asymptotic normality. The corresponding theoretical risk is $R(\theta)=0.5 v^{-1} (\E Z-\phi)^2$, minimized at $\phi=\E Z$.

The model here is partially identified, since the quasi-likelihood $e^{-\lambda R_n(\theta)}$ only depends on the reduced-form parameter $\phi$. The data can only identify $\phi$. Given $\phi$, the structural parameter $\omega$ can still be anywhere from the prior support of $\pi(\omega|\phi) \propto \sum_j \pi(M_j)\pi(\omega|M_j)\pi(\phi|\omega)$, which is supported on $\Omega(\phi)=[\phi,\phi+1]\cap[0,4]$. Here $\Omega(\phi)$ is called the {\it identification region} for $\omega$ given $\phi$. This is related to the minimizer of of the theoretical risk of $R$ when $R$ is regarded as a function of $\theta=(\omega,\phi)$, even if it depends really on $\phi$ only. Suppose $R$ has a unique minimizer $\phi=\phi^*$ (the ``true" $\phi$), then attaching all possible $\omega$ values in $\Omega(\phi^*)$,
we have
$$\arg\min_\theta R(\theta) = \{\phi^*\}\times \Omega(\phi^*).$$
Suppose the true $\phi^*=3.6$. Then the identification region for $\omega$ is $\Omega(3.6)=[3.6,4]$, and $\arg\min_\theta R(\theta) =\{3.6\}\times [3.6,4]$.

Model $M_1$ is ``incompatible'' with data, in the sense that its prior cannot reach the minimum theoretical risk for $R(\theta)$. The proposed prior on $\omega$ does not allow $\phi=\phi^*$, the risk minimizer and the true $\phi$. In other words, the prior support of $\pi(\theta|M_1)= \pi(\phi|\omega)\pi(\omega|M_1)$ does not intersect $\arg\min_{\theta} R(\theta) =\{\phi^*\}\times \Omega(\phi^*)$, since the support of $\pi(\omega|M_1)$ is $\{3\}$, which does not intersect with $\Omega(\phi^*)=[3.6,4]$.

Model $M_2$ is ``compatible'' with data, in the sense that its prior can  reach the minimum theoretical risk $R$. The proposed prior on $\omega$ does allow $\phi=\phi^*$, the risk minimizer and the true $\phi$.) In other words, the prior support of $\pi(\theta|M_2)= \pi (\phi|\omega)\pi(\omega|M_2)$ intersects $\arg\min_\theta R(\theta) = \{\phi^*\}\times \Omega(\phi^*)$, since the support of $\pi(\omega|M_2)$ is $[0,4]$, which intersects with $\Omega(\phi^*)=[3.6,4]$.

This simple example will be generalized in the next section \ref{oraclepid}, where there can be more than two  model candidates and the quasi-posterior can also involve more than two parameters. We hope that with Bayesian model averaging, incompatible models can have small posterior probability asymptotically, so that the posterior from model averaging will be as good as the oracle posterior, which assumes that we knew beforehand and had only used those models that are compatible with data.

\subsection{Bayesian model averaging and oracle properties with partial identification}\label{oraclepid}

We first derive oracle properties for model selection and BMA in the general framework of quasi-posterior as defined in \eqref{quasipost0}. Later we will consider the special case of partial identification described in \citet{ms12}.


Define the index set $\Jt=\left\{j\geq 1: \inf_{\theta\in\Theta_j}  R(\theta) = \inf_{\theta\in\Theta}  R(\theta)\right\}$, which includes all model indexes under which the global minimum risk can be reached.  These models will be called ``compatible models''. With partial identification, it is important to allow all compatible models in consideration, and not to exclusively favor one compatible model, even if it is the simplest model with the lowest model complexity. An alternative approach could use a dimensional penalty to favor the simplest compatible model, but this could miss true values of the parameter $\theta$ due to partial identification, as discussed in an earlier technical report \citet{jl15} Section 6.6.2.
Another example that illustrates this kind of subtlety is described as a technical detail in a supplementary material of the current paper.

In response to this subtlety with partial identification, we will group all the compatible models together to form our ``true" model $M^*=\left\{M_j: j\in \Jt\right\}$. Then $\pi(\theta, M^*|\DD) \propto e^{-\lambda R_n(\theta) } \sum_{j\in \Jt} \pi(\theta, M_j)$. The resulting joint prior on $\theta$ and $M_j$ can be rewritten as $\pi(\theta, M^*) = \sum_{j\in \Jt} \pi(\theta, M_j) = \pi(M^*)\pi(\theta|M^*)$, where $\pi(M^*)=\sum_{j\in\Jt} \pi(M_j)$, and $\pi(\theta|M^*) = \sum_{j\in \Jt} \pi(\theta|M_j) \cdot \pi(M_j) / \sum_{j\in \Jt} \pi(M_j)$ is a mixture prior for $\theta$ conditional on the composite true model $M^*$.

All incompatible models are indexed by $j\in \Jw$.  For incompatible models, we assume the quantity $\gamma = \inf_{j\in \Jw}\inf_{\theta\in\Theta_j} R(\theta)- \inf_{\theta\in\Theta}  R(\theta)$ to be a positive constant, which holds true if there is a fixed number of candidate models. This $\gamma$ is exactly the same $\gamma$ used in Proposition \ref{gibbs1}.  We can derive an upper bound for the posterior mis-selection probability $1-\pi(M^*|\DD)$ (where $M^*=\left\{M_j: j\in \Jt\right\}$) as exponentially small in $\lambda$ from Proposition \ref{gibbs1}, which leads to the Bayesian oracle properties \ref{O1} and \ref{O2} in Section 2. The oracle posterior here is still $\pi(\theta|M^*,\DD)$, conditional on compatible models only.

We make the following assumptions.

\begin{itemize}
\item[(A1)] $\lambda \succ 1$ as $n\rightarrow\infty$.
\item[(A2)] $\sup_{\theta\in \Theta} |R_n(\theta)- R(\theta)|=o_p(1)$ as $n\to\infty$.
\item[(A3)] $\Pi\left(\left\{\theta: R(\theta)-\inf_{\theta\in\Theta} R(\theta) < a \right\}\right)>0$ for any small $a>0$.
\item[(A4)] $\gamma = \inf_{j\in \Jw}\inf_{\theta\in\Theta_j} R(\theta)- \inf_{\theta\in\Theta} R(\theta)$ is a positive constant.
\end{itemize}

Assumption (A1) is true when $\lambda\propto n$. When $R_n$ is a sample average of independently and identically distributed data, we can take the theoretical risk in (A2) to be the expectation $R(\theta)=\E R_n(\theta)$ over the true distribution of the randomly generated data. Then (A2) can be satisfied due to a uniform law of large numbers, which holds, e.g., when the entire parameter space $\Theta$ is compact and the risk functions are stochastically equicontinuous (see, e.g, \citealt{n91}). \footnote{Assumption (A2) may also be satisfied when $R_n$ is not an average itself, but is a function of some sample averages, such as is easy too verify for the example in Section \ref{simpleexample}. In fact it is easy to check that all conditions are valid for that example assuming that ${\rm var}(Z)>0$.}  Assumption (A4) is true when the number of candidate models is fixed. Regarding (A3), suppose the prior support $\Theta$ is compact and contains a risk minimizer of $R$ in its interior. Then a small enough neighborhood of this risk minimizer will have positive prior $\pi$ and can have risk $R(\theta)$ being arbitrarily close to the minimum risk, if $R(\theta)$ is continuous in $\theta$.

We can summarize the analysis above formally in the following theorem.
\begin{thm} \label{partselect}
Assume that (A1)-(A4) hold and $M^*=\left\{M_j: j\in \Jt\right\}$. Then the total variation distance between the distributions $\pi(\theta|\DD)$ and $\pi(\theta|M^*,\DD)$ is $o_p(1)$ as the sample size $n\to\infty$, i.e., the global model selection consistency (Property \ref{O1}) and the Bayesian model averaging oracle property (Property \ref{O2}) both hold.
\end{thm}

The proof of Theorem \ref{partselect} shows $\pi(\Jw)=o_p(1)$ by applying Proposition \ref{gibbs1}. Therefore even though it is impossible to point identify the minimizer of the theoretical risk, we can still have a similar form of Bayesian oracle properties by selecting all the compatible models. As a result, the posterior inference based on model averaging is asymptotically equivalent to the posterior inference based on only those compatible models weighted by their priors.

\newtheorem{rem}{Remark}
\newtheorem{cor}{Corollary}

The above Theorem \ref{partselect} is very general. \citet{ms12} considered a special case where $-\lambda R_n(\theta)$ is  the log likelihood function. Also, the ``combined'' parameter can be decomposed as $\theta=(\omega,\phi)$, where $\omega$ is a structural parameter of interest and $\phi$ is a reduced-form parameter that is identified by data. The candidate models impose different priors on the structural parameter $\omega$, so that $\pi(\theta,M_j)=\pi(M_j)\pi(\theta|M_j)=\pi(M_j)\pi(\omega|M_j)\pi(\phi|\omega)$. A simple example of this kind of parametrization and the corresponding prior distribution is described in Section \ref{simpleexample}, using a quasi-likelihood derived from asymptotic normality. For such situations when only the structural parameter $\omega$ is of primary interest, the BMA oracle property \ref{O2} for the marginal posterior on $\omega$ also holds:
\begin{cor}\label{marginal}
Under the assumptions made for Theorem \ref{partselect}, the BMA oracle property \ref{O2} holds marginally for the structural parameter of interest $\omega$, i.e., $\int|\pi(\omega|\DD)-\pi(\omega|M^*,\DD)|d\omega=o_p(1)$, if $\omega$ is a sub-vector of the combined parameter $\theta$.
\end{cor}

So far we have discussed Property \ref{O1} (for global model selection consistency) and Property \ref{O2} (for the oracle property with BMA). There is an important exception here: Property \ref{O3} for MAP model selection is not guaranteed in this partially identified model. This is because here the true model $M^*$ is effectively the set of all compatible models which is possibly a nonsingleton, and the proof of Proposition \ref{tv3} does not go through. When there are two or more compatible models, the MAP model selection may only choose one compatible model and neglect all the other ones. Posterior inference based on the MAP model may be different from using the oracle posterior given all the compatible models and may end up missing the true value of a point parameter. We will describe this as a technical detail with a simple example in a supplementary material.

Regarding the mean oracle property \ref{O4}, we conjecture that it usually holds for the structural parameter of interest, as will be discussed as some additional technical details in the supplementary material.

\section{Discussion}
We have established a fundamental relation between three different topics:
Bayesian model averaging, model selection consistency, and oracle performance in posterior distribution.
The relatively basic property of model selection consistency is shown to imply a seemingly more advanced distributional result, the oracle property. The result is very simple and general. Unlike some previous Bayesian oracle properties discussed in special cases such as \citet{ir11}, and \citet*{csv15}, who consider linear models, and \citet{hp12} and \citet{lj16} who consider identifiable models with standard limiting distributions, the current work is completely free from any restriction on the type of prior or (quasi-)likelihood function used, or even from any restriction on the limiting distribution of the oracle posterior.

For applications, we considered two classes of models with nonstandard limiting distributions studied in \citet{ms12} and \citet*{jpw15}. They involve partial identifiability or nonstandard rates of convergence, but we can still show the Bayesian oracle properties, which suggest that Bayesian model averaging  can be applied to their methods and work well for Bayesian inference of the unknown point parameter.  On the other hand, we suspect that  model selection based on MAP may not be reliable  for the partial identification example and may miss reasonable models (see a discussion after Corollary \ref{marginal}).

When the model is misspecified, the model that minimizes the theoretical risk $R$ plays the role of the true model in our theory. Our oracle property will imply that the quasi-posterior based on BMA will converge to the quasi-posterior based on the minimum risk model, asymptotically. \citet{gvo14} discovered suboptimal predictive performance when a homoscedastic linear model is misspecified. Their numerical experiments seem to indicate that the performance of BMA still converges to the performance of the true model eventually, albeit with a much larger sample size compared to the correctly specified case. This indicates a much slower convergence speed of BMA when the models are misspecified. Our current paper only addresses the limiting distributional behavior of BMA and BMS, but not their convergence speed. As a possible future work, we may consider extending our theory in Section 3.1 to study the convergence speed in the presence of model misspecification and how the convergence depends on the temperature parameter, as discussed in \citet{gvo14}.

Given the success of the frequentist oracle properties studied by \citet{fl01}, we expect that the Bayesian version should also have applications in a wide variety of situations, in addition to the examples discussed in this paper. For example, the relationships described in Section 2 and 3 can be generalized to models with increasing or high dimensions, and potentially to other nonstandard model selection problems with appropriate conditions on the priors.  For instance, \citet{dp17} have developed a generalized version of BIC type approximation for the class of singular models (such as factor models), where the posterior model probability does not allow a quadratic approximation and results in an extra $\ln\ln n$ term in the BIC approximation \eqref{bic11}. Our Bayesian oracle properties may also apply to these singular models. In addition, in the context with partial identification, our paper only considered inference about the partially identified point parameter, following the approach of, e.g., \cite{p98}, \citet{ms12}, and \citet{gus15}.  It may also be of interest to consider inference about the fully identified ``set parameter'', following the approach of, e.g., \citet{w13}, \citet{kt16}, and \citet*{cct16}, and develop similar oracle properties for Bayesian model selection or model averaging.


\section{Technical Details}
\subsection{Proof of propositions}\label{prfs}

\begin{proof}[Proof of Proposition \ref{tvmisprob}]

For any event $A$ and $B$, we have
\begin{align*}
&\Pi(A|\DD)=\Pi(A|B,\DD)\Pi(B|\DD)+\Pi(A|B^c,\DD)\Pi(B^c|\DD), \\
&\Pi(A|B,\DD)=\Pi(A|B,\DD)\Pi(B|\DD)+\Pi(A|B,\DD)\Pi(B^c|\DD).
\end{align*}
Therefore
$$|\Pi(A|\DD)-\Pi(A|B,\DD)| = |\Pi(A|B^c,\DD)-\Pi(A|B,\DD)|\Pi(B^c|\DD)\leq \Pi(B^c|\DD)$$
for any $A$. Taking supremum over all event $A$ and setting event $B=\{M=M^*\}$ lead to the proof. (Note that \citealt*{csv15} used a double-sized upper bound in proving their Theorem 6 in the context of Bayesian linear regression.)
\end{proof}

\begin{proof}[Proof of Proposition \ref{tv3}]

The MAP choice $\widehat M$ satisfies $\pi(\widehat M|\DD)\geq \pi(M^*|\DD)$ by definition. In the proof of Proposition \ref{tvmisprob} above, we can replace $M^*$ by $\widehat M$ and obtain that $\sup_A|\Pi(A|\DD)-\Pi(A|\widehat M,\DD)|\leq 1-\pi(\widehat M|\DD)$.
The right hand side is at most $1-\pi(M^*|\DD)$ since $\pi(\widehat M|\DD)\geq \pi(M^*|\DD)$. Now combining this with the result of Proposition \ref{tvmisprob} using the triangle inequality leads to the conclusion.
\end{proof}

\begin{proof}[Proof of Proposition \ref{bicproperty}]

Let $p(M_j)=\int_{\Theta_j} e^{-\lambda R_n(\theta)} d\pi(\theta|M_j)$. Under the assumption (ii), we have that for any minimum-risk model $M_j\neq M^*$,
\begin{align}\label{bic12}
&-\ln p(M_j)=  \lambda R_n(\theta^*) + \frac{d_j\ln\lambda}{2} + O_p(1), \nonumber \\
&-\ln p(M^*) = \lambda R_n(\theta^*) + \frac{d^*\ln\lambda}{2} + O_p(1).
\end{align}
Taking the difference between these two equations gives
\begin{align*}
-\ln p(M_j)/p(M^*) = \frac{(d_j-d^*)\ln \lambda }{2} +O_p(1).
\end{align*}
Due to the assumptions (iii) and (v), $O_p(1)$ is negligible compared to the $\ln\lambda$ term. Therefore, for any minimum-risk model $M_j$, there exists a constant $C_1>0$ such that
\begin{align}\label{mtrue1}
p(M_j)/p(M^*)\leq C_1\lambda^{-(d_j-d^*)/4}\leq C_1\lambda^{-1/4}.
\end{align}
We notice that from \eqref{bic12}, $p(M^*)=\exp\left[-\lambda R_n(\theta) - d^*\ln \lambda/2 +O_p(1)\right]$. For any non-minimum-risk model $M_j$, we can use the assumption (iv) to obtain that for some constant $C_2>0$,
\begin{align}\label{mwrong1}
p(M_j) &= \int_{\Theta_j} e^{-\lambda S_n(\theta) - \lambda[R(\theta)-R(\theta^*)] - \lambda R_n(\theta^*)} d\pi(\theta|M_j) \nonumber \\
&\leq \int_{\Theta_j} e^{-\lambda S_n(\theta) } d\pi(\theta|M_j)\cdot e^{-\lambda\gamma_j} \cdot p(M^*) e^{\frac{d^*\ln \lambda}{2}+O_p(1)} \nonumber  \\
&\leq C_2 \lambda^{d^*/2}e^{-\lambda \gamma_j} p(M^*).
\end{align}
Since $\gamma_j\succeq 1$ and $\max_{j\geq 1} d_j$ is upper bounded by constant, the exponential rate $e^{-\lambda \gamma_j}$ dominates the polynomial rate $\lambda^{d^*/2}$. Furthermore, from the assumption (i), we also have that the prior ratio $\pi(M_j)/\pi(M^*)$ is lower and upper bounded by constants for any model $M_j$. Therefore, from \eqref{mtrue1} and \eqref{mwrong1}, we have that
\begin{align*}
& 1-\pi(M^*|\DD) \\
& = \frac{\sum_{M_j\neq M^*} \pi(M_j) p(M_j)}{\sum_{M_j\neq M^*} \pi(M_j) p(M_j) + \pi(M^*) p(M^*)} = \frac{\sum_{M_j\neq M^*} \frac{\pi(M_j)}{\pi(M^*)} \frac{p(M_j)}{p(M^*)}}{\sum_{M_j\neq M^*} \frac{\pi(M)}{\pi(M^*)} \frac{p(M)}{p(M^*)} + 1} \\
&\leq 1- \Big[\sum_{M_j\neq M^*\text{ and }M_j\text{ is minimum-risk}} \frac{\pi(M_j)}{\pi(M^*)}C_1\lambda^{-1/4} \\
&~~~ + \sum_{M_j \text{ is non-minimum-risk}} \frac{\pi(M_j)}{\pi(M^*)} C_2 \lambda^{d^*/2}e^{-\lambda \gamma_j} +1 \Big]^{-1} \\
& =o_p(1).
\end{align*}
Therefore the global model selection consistency (Property \ref{O1}) is proved. By Propositions \ref{tvmisprob} and \ref{tv3}, the Bayesian oracle properties for BMA (Property \ref{O2}) and BMS (Property \ref{O3}) also hold.

Furthermore, if the assumption (vi) holds, then $\E(\theta|M_j,\DD)$ exists and its $L_2$ norm is uniformly bounded by some constant $C_3>0$ for all models $M_j$ since $\Theta$ is compact. Based on the assumption (vi), for any true model $M_j$, a triangle inequality yields
\begin{align*}
\left\|\E(\theta|M_j,\DD)- \E(\theta|M^*,\DD)\right\| &\leq \left\|\E(\theta|M_j,\DD)-\theta^*\right\| +\left\|\E(\theta|M^*,\DD)-\theta^*\right\|\\
& = \left[1+O_p(1)\right] \left\|\E(\theta|M^*,\DD)-\theta^*\right\|.
\end{align*}
This together with \eqref{mtrue1} and \eqref{mwrong1} implies that
\begin{align}\label{meanoraclebound}
& \left\|\E(\theta|\DD)-\E(\theta|M^*,\DD) \right\| = \left\|\sum_{j\geq 1, M_j\neq M^*}\pi(M_j|\DD)\left[\E(\theta|M_j,\DD)-\E(\theta|M^*,\DD)\right]\right\| \nonumber \\
& \leq \left\|\sum_{M_j\neq M^*\text{ and }M_j\text{ is minimum-risk}}\pi(M_j|\DD)\left[\E(\theta|M_j,\DD)-\E(\theta|M^*,\DD)\right]\right\| \nonumber \\
& ~~ + \left\|\sum_{M_j \text{ is non-minimum-risk}}\pi(M_j|\DD)\left[\E(\theta|M_j,\DD)-\E(\theta|M^*,\DD)\right]\right\| \nonumber \\
& \leq \left[1+O_p(1)\right] \left\|\E(\theta|M^*,\DD)-\theta^*\right\|\cdot \sum_{M_j\neq M^*\text{ and }M_j\text{ is minimum-risk}} \pi(M_j|\DD) \nonumber  \\
& ~~ + 2C_3\sum_{M_j \text{ is non-minimum-risk}}\pi(M_j|\DD) \nonumber \\
& \leq \left[1+O_p(1)\right] \left\|\E(\theta|M^*,\DD)-\theta^*\right\| \cdot \sum_{M_j\neq M^*\text{ and }M_j\text{ is minimum-risk}} \frac{\pi(M^*|\DD)\pi(M_j)}{\pi(M^*)}C_1\lambda^{-1/4} \nonumber \\
& ~~ + 2C_3\sum_{M_j \text{ is non-minimum-risk}} \frac{\pi(M^*|\DD)\pi(M_j)}{\pi(M^*)} C_2 \lambda^{d^*/2}e^{-\lambda \gamma_j}.
\end{align}
First notice that the number of models is bounded and all prior ratios $\pi(M_j)/\pi(M^*)$ are upper bounded by constants, according to the assumption (i). Since $\lambda \to\infty$ by the assumption (v), and $\pi(M^*|\DD)$ is in $[0,1]$, we have that the first term in \eqref{meanoraclebound} is $o_p(1)\cdot\left\|\E(\theta|M^*,\DD)-\theta^*\right\|$ as $n\to \infty$. Due to the assumptions (iv) and (vi), $\lambda^{d^*/2}e^{-\lambda \gamma_j}$ decays exponentially fast in $n$, while $\left\|\E(\theta|M^*,\DD)-\theta^*\right\|$ decays polynomially in $n$ with a rate no faster than $\epsilon_n$. Hence it is clear that the second term in \eqref{meanoraclebound} is also of order $o_p(1)\cdot\left\|\E(\theta|M^*,\DD)-\theta^*\right\|$. Therefore, both terms in \eqref{meanoraclebound} are $o_p(1)\cdot\left\|\E(\theta|M^*,\DD)-\theta^*\right\|$, and we have proved the mean oracle property (Property \ref{O4}).
\end{proof}

\begin{proof}[Proof of Proposition \ref{gibbs1}]

First it is clear that on all models not equal to $M^*$, $R(\theta) - \inf_{\theta\in \Theta} R(\theta) \geq \inf_{\theta\in \Theta, M\neq M^*} R(\theta) - \inf_{\theta\in \Theta} R(\theta) =\gamma$. Hence
\begin{equation}\label{msrisk}
1-\pi(M^*|\DD)\leq \Pi\left(\left\{\theta: R(\theta) \geq \inf_{\theta\in \Theta} R(\theta) + \gamma\right\}\Big|\DD\right).
\end{equation}

We then show that for any bounded measurable function $h(\theta)$,
\begin{equation}\label{gibbslim}
\ln {\E}[h(\theta)|\DD] \leq \frac{1}{2}\ln {\E}_{\infty}[h^2(\theta)] -\lambda u,
\end{equation}
where ${\E}[h(\theta)|\DD]=\int_{\theta\in \Theta}h(\theta) \pi(\theta|\DD)d\theta$, ${\E}_{\infty}[h^2(\theta)]=\int_{\theta\in \Theta} h^2(\theta) \pi_{\infty}(\theta)d\theta$, $u$ is defined as in Proposition \ref{gibbs1}.
To see why \eqref{gibbslim} holds true, we recall the definitions of a quasi-posterior $\pi(\theta|\DD)$ and its ``limiting posterior" $\pi_{\infty}(\theta)$:
\begin{align*}
& {\E}[h(\theta)|\DD] =\frac{\int_{\theta\in \Theta} e^{-\lambda R_n(\theta)} h(\theta) \pi(\theta) d\theta} {\int_{\theta\in \Theta}  e^{-\lambda R_n(\theta)} \pi(\theta) d\theta} \\
& =\frac{\int_{\theta\in \Theta} e^{-\lambda [R_n(\theta)-R(\theta)]}h(\theta) \pi(\theta) d\theta} {\int_{\theta\in \Theta} e^{-\lambda [R_n(\theta)-R(\theta)]} \pi(\theta) d\theta}.
\end{align*}
Then we apply the Jensen's inequality to the denominator and apply the Cauchy-Schwarz inequality to the numerator to obtain that
\begin{align*}
& {\E}[h(\theta)|\DD] \leq \frac{\sqrt{\int_{\theta\in \Theta} e^{-2\lambda[R_n(\theta)-R(\theta)]}\pi_{\infty}(\theta) d\theta} \sqrt{\int_{\theta\in \Theta} h^2(\theta) \pi_{\infty}(\theta)d\theta}} {e^{-\lambda \int_{\theta\in \Theta} [R_n(\theta)-R(\theta)] \pi_{\infty}(\theta)d\theta}} \\
& =\sqrt{\int_{\theta\in \Theta} e^{-2\lambda \left[(R_n(\theta)-R(\theta)) -\int_{\theta\in \Theta}(R_n(\theta)-R(\theta))\pi_{\infty}(\theta)d\theta\right] }\pi_{\infty}(\theta)d\theta} \sqrt{\int_{\theta\in \Theta} h^2(\theta) \pi_{\infty}(\theta)d\theta},
\end{align*}
which leads to \eqref{gibbslim}. Then we take $h=I(A)$ for a measurable set $A$ and obtain that
\begin{equation}\label{gibbslimevent}
\ln \Pi(A|\DD) \leq \frac{1}{2} \ln \Pi_{\infty}(A) -\lambda u.
\end{equation}
Set $A=\{\theta: R(\theta)-\inf_{\theta\in \Theta} R(\theta)\geq \gamma\}$ and use the definition of $r$ in Proposition \ref{gibbs1}:
\begin{align*}
& \Pi_{\infty}(A) = \frac{\int_{\theta\in \Theta} e^{-\lambda [R(\theta)-\inf_{\theta\in \Theta} R(\theta)] } I(A) \pi(\theta)d\theta}{\int_{\theta\in \Theta} e^{-\lambda [R(\theta)-\inf_{\theta\in \Theta} R(\theta)] }\pi(\theta)d\theta} \\
& =\int_{\theta\in \Theta} e^{-\lambda [R(\theta)-\inf_{\theta\in \Theta} R(\theta)-r] } I(A) \pi(\theta)d\theta \leq e^{-\lambda (\gamma-r)}.
\end{align*}
Then applying this upper bound of $\Pi_{\infty}(A)$ to \eqref{gibbslimevent} and using \eqref{msrisk} leads to the proof.
\end{proof}

\subsection{Proof of theorems}

\begin{proof}[Proof of Theorem \ref{cubicroot}]

This is similar to the proof of the proposition in Section \ref{bic01}.
We proceed in two steps: first we show the global model selection consistency (Property \ref{O1}) (which implies Properties \ref{O2} and \ref{O3} by Propositions \ref{tvmisprob} and \ref{tv3}), and then we show the mean oracle property (Property \ref{O4}).
\vspace{.3cm}

\noindent Step 1: Show the global model selection consistency (Property \ref{O1}).
\vspace{.2cm}

With the model selection prior, the quasi-posterior in \eqref{quasipost0} can be rewritten as
\begin{align}\label{cpost1}
\pi(\theta|M_j,\DD) & = \frac{\pi(\theta|M_j,\DD) \exp\{-\lambda R_n(\theta)\} }
{ \pi(M_j|\DD)}, \nonumber \\
\pi(M_j|\DD) & = \frac{\pi(M_j)p(M_j)}{\sum_{l\geq 1} \pi(M_l) p(M_l)}, \nonumber \\
p(M_j) & = \int_{\Theta_j} \pi(\theta|M_j) \exp\{-\lambda R_n(\theta)\} d\theta,
\end{align}
for any model $M_j$ as a coordinate subspace of $\Theta\cap \mathbb{R}^p$.
\vspace{.2cm}

We group all the models $M_j$ that are different from the true model $M^*$ into 2 separate groups. Group 1 contains the models that include $M^*$ as a strict submodel, i.e. $M_j\supset M^*$ and $d_j\geq p^*+1$. (Group 1 does not exist if $p^*=p$). Group 2 contains the models that miss at least one component of $M^*$, i.e. $M^*\backslash M_j \neq \emptyset$.

Define the localize parameter $t=\sqrt{\lambda}(\theta-\theta^*)$ for all $\theta\in \Theta_j$. Define the quantities
\begin{align}\label{SRdef}
\widetilde S_n(t) &= n^{1/2}\lambda^{1/4} \left\{\left[R_n\left(\theta^*+\frac{t}{\sqrt{\lambda}}\right) -  R\left(\theta^*+\frac{t}{\sqrt{\lambda}}\right)\right] - \left[R_n(\theta^*)-R(\theta^*)\right]\right\}, \nonumber \\
\widetilde R(t) &= \lambda \left[R\left(\theta^*+\frac{t}{\sqrt{\lambda}}\right)-R(\theta^*)\right],
\end{align}
where $\widetilde S_n(t)$ and $\widetilde R(t)$ are the rescaled and centered versions of $S_n(\theta)$ (defined in Section \ref{bic01}) and $R(\theta)$.

Then the posterior model probability $p(M)$ has the expression
\begin{align}\label{pmt2}
p(M_j) & = \lambda^{-d_j/2} e^{-\lambda R_n(\theta^*)} \int_{\mathcal{T}_{nj}} \pi_n(t|M_j)
\exp\left\{-\left[n^{-1/2}\lambda^{3/4}\widetilde S_n(t) + \widetilde R(t) \right]\right\} dt,
\end{align}
where $\mathcal{T}_{nj} = \left\{t=\sqrt{\lambda} (\theta-\theta^*): \theta\in \Theta_j\right\}$ and
$\pi_n(t|M_j) = \pi(\theta^*+t/\sqrt{\lambda}~ |M_j)$.

For the true model $M^*$ and any model $M_j$ in Group 1, from the proof of Theorem 1 part (iii) of \citet*{jpw15}, we have
\begin{align}\label{pmg1}
p(M_j) & = \lambda^{-d_j/2}e^{-\lambda R_n(\theta^*)} (2\pi)^{d_j/2} \left|\dett(V_{M_j})\right|^{-1/2} \pi(\theta^*|M_j)(1+o_p(1)),
\end{align}
where $\pi=3.1415926\ldots$, $V_{M_j}$ is the principle minor of $V$ in (C3) restricted to the components in $M$, and $\dett(\cdot)$ represents the determinant. According to (C3), since the eigenvalues of $V$ are bounded from below and above by positive constants, so are the eigenvalues of $V_{M_j}$. According to (C5), $\pi(\theta^*|M_j)$ are bounded from above and $\pi(\theta^*|M^*)$ is bounded from below by a positive constant. Therefore for any model $M_j$ in Group 1, for all sufficiently large $n$, with probability approaching 1,
\begin{align}\label{ratiog1}
\frac{p(M_j)}{p(M^*)} & =  (1+o_p(1)) (2\pi/\lambda)^{(d_j-p^*)/2} \left[\frac{\left|\dett(V_{M^*})\right|}{\left|\dett(V_{M_j})\right|}\right]^{1/2} \cdot \frac{\pi(\theta^*|M_j)}{\pi(\theta^*|M^*)} \nonumber \\
&\leq c_1 \lambda^{-(d_j-p^*)/2}\leq c_1 \lambda^{-1/2},
\end{align}
where $c_1>0$ is a constant, and $d_j-p^*\geq 1$ for any model $M_j$ in Group 1.
\vspace{.2cm}

For any model $M_j$ in Group 2, since $M_j$ misses at least one nonzero component of $\theta^*$, by (C1) we have $\|\theta-\theta^*\|\geq c_{\theta}$ for all $\theta\in \Theta_j$. By the compactness of $\Theta$ from (C1), the uniqueness of $\theta^*$ from (C2), and the continuity of $R(\theta)$ from (C3), there exists a constant $\gamma>0$ such that $R(\theta)-R(\theta^*)\geq \gamma$ whenever $\|\theta-\theta^*\|\geq c_{\theta}$. Therefore, for any model $M_j$ in Group 2, $\inf_{\Theta_j} R(\theta) - R(\theta^*) \geq \gamma$, and for $\widetilde R(t)$ defined in \eqref{SRdef}, we have $\widetilde R(t)\geq \lambda \gamma$ for all $t\in \mathcal{T}_{nj}$.

Under (C4), (C6), (C7) and (C8), Lemma B.2 of \citet*{jpw15} shows that
\begin{align} \label{gplimit}
& \widetilde S_n(t) \rightarrow_{L^{\infty}_L(\mathbb{R}^p)} \mathbb{G}(t),
\end{align}
where $\mathbb{G}(t)$ is a Gaussian process indexed by $t$ with zero mean and covariance kernel $H(s,t)$ as defined in (C4), and $L^{\infty}_L(\mathbb{R}^p)$ is the space of all locally bounded functions on compacta. Based on this limit and (C4), their Lemma B.6 furthers shows that
\begin{align}\label{supsnt}
\sup_{t\in \mathbb{R}^p} \left[\left|\widetilde S_n(t)\right|-c\|t\|\right] = O_p(1),
\end{align}
where $c>0$ is an arbitrary constant.

Now we can bound $p(M_j)$ in \eqref{pmt2} for $M_j$ in Group 2, in a similar manner to the proof of Proposition \ref{bicproperty}. Let $c_{\pi}>0$ be a constant upper bound for $\pi(\theta|M_j)$ by (C5). Then for any fixed $c>0$ we have
\begin{align*}
& p(M_j) \leq \lambda^{-\frac{d_j}{2}}e^{-\lambda R_n(\theta^*)} \int_{\mathcal{T}_{nj}} c_{\pi}
\exp\left\{-\left[n^{-1/2}\lambda^{3/4}\widetilde S_n(t) + \widetilde R(t) \right]\right\} dt \\
&\leq c_{\pi}  \lambda^{-\frac{d_j}{2}}  e^{-\lambda R_n(\theta^*)- \lambda\gamma} \int_{\mathcal{T}_{nj}}
\exp\left(n^{-1/2}\lambda^{3/4}\left\{\left|\sup_{t\in \mathbb{R}^p} \left[\left|\widetilde S_n(t)\right|-c\|t\|\right] \right|+ c\|t\|\right\}\right) dt.
\end{align*}
In the exponent of the integrand, $n^{-1/2}\lambda^{3/4} = o(1)$, $\|t\|\leq 2\sqrt{\lambda} \sup_{\theta\in \Theta}\|\theta\| \prec \lambda$ due to the compactness of $\Theta$ in (C1), and the supremum term is $O_p(1)$ by \eqref{supsnt}. Thus for any model $M$ in Group 2, for all sufficiently large $n$, the following inequality holds with probability approaching 1:
\begin{align}\label{pmg2}
p(M_j) & \leq c_{\pi}\lambda^{-\frac{d_j}{2}} e^{-\lambda R_n(\theta^*)}\cdot \voll(\mathcal{T}_{nj}) \exp\left(-\frac{3}{4}\lambda\gamma \right) \nonumber \\
&\leq \voll(\Theta) c_{\pi} \exp\left(-\lambda R_n(\theta^*)-\frac{3}{4}\lambda\gamma \right),
\end{align}
where $\voll(\cdot)$ represents the volume of a set. Therefore by \eqref{pmg1} and \eqref{pmg2}, for any model $M_j$ in Group 2, for all sufficiently large $n$, the following bound on the posterior odds ratio holds true with probability approaching 1:
\begin{align}\label{ratiog2}
\frac{p(M_j)}{p(M^*)} & \leq \frac{\voll(\Theta) c_{\pi} \exp\left(-\lambda R_n(\theta^*)-\frac{3}{4}\lambda\gamma \right)} {\lambda^{-\frac{p^*}{2}}e^{-\lambda R_n(\theta^*)} (2\pi)^{p^*/2} \left|\dett(V_{M^*})\right|^{-1/2} \pi(\theta^*|M^*)} \nonumber \\
&\leq c_2 \lambda^{\frac{p}{2}}  \exp\left(-\frac{3}{4}\lambda\gamma \right) \leq c_2 \exp\left(-\frac{1}{2}\lambda\gamma \right)
\end{align}
where $c_2>0$ is a constant. The last upper bound goes to zero as $n\to \infty$ because $n^{2/5}\prec \lambda \prec n^{2/3}$.

We combine \eqref{cpost1}, \eqref{ratiog1} and \eqref{ratiog2} together and obtain that
\begin{align*}
& 1- \pi(M^*|\DD) \\
& = \frac{\sum_{M_j \neq M^*} \pi(M_j) p(M_j)}{\sum_{M_j\neq M^*} \pi(M_j) p(M_j) + \pi(M^*) p(M^*)} = \frac{\sum_{M_j\neq M^*} \frac{\pi(M_j)}{\pi(M^*)} \frac{p(M_j)}{p(M^*)}}{\sum_{M_j\neq M^*} \frac{\pi(M_j)}{\pi(M^*)} \frac{p(M_j)}{p(M^*)} + 1} \\
&\leq 1- \left[\sum_{M_j \in \text{Group 1}} \frac{\pi(M_j)}{\pi(M^*)}c_1\lambda^{-1/2} + \sum_{M_j \in \text{Group 2}} \frac{\pi(M_j)}{\pi(M^*)} c_2 \exp\left(-\frac{1}{2}\lambda\gamma \right) +1 \right]^{-1} \\
& =o_p(1)
\end{align*}
as $n\to \infty$, since the prior ratios $\pi(M_j)/\pi(M^*)$ are uniformly bounded from above by constant from (C5). This has proved the global model selection consistency (Property \ref{O1}).
\vspace{.4cm}

\noindent Step 2: Show the mean oracle property (Property \ref{O4}).
\vspace{.2cm}

We use the relation \eqref{detailrate} to show the mean oracle property \ref{O4}. For the true model $M^*$ and any model $M_j$ in Group 1, the proof of Theorem 1 part (iii) in \citet*{jpw15} has shown that
\begin{align}\label{mdg1}
n^{1/2}\lambda^{-1/4}\left[\E(\theta|M_j,\DD)-\theta^*\right] & =
\int_{\mathcal{T}_{nj}} \frac{\pi_n(t|M_j)}{\pi(\theta^*|M_j)} \cdot t \widetilde S_n(t) \phi_{V_{M_j}}(t) dt + o_p(1),
\end{align}
where $\phi_{V_{M_j}}(t)$ is the density of $N(0,V_{M_j}^{-1})$. Because \eqref{gplimit} shows that $\widetilde S_n(t)$ converges to a Gaussian process with covariance kernel $H(t,s)$, the first term of integral on the right hand side of \eqref{mdg1} converges in distribution to $N(0,\mathcal{V}_{M_j})$, where $\mathcal{V}_{M_j} =\qquad$ $\iint_{\mathcal{T}_{nj}^2} ts^{\top}H(t,s)\phi_{V_{M_j}}(t)\phi_{V_{M_j}}(s)dtds$. Although Theorem 1 part (iii) of \citet*{jpw15} has an extra asymptotic bias term in the normal limit, it follows from their remarks after their Theorem 1 that under $n^{2/5}\prec \lambda \prec n^{2/3}$ and our conditions about bounded derivatives in (C3) and (C5), the asymptotic bias term is also a negligible $o_p(1)$, which can be absorbed into the $o_p(1)$ on the right hand side of \eqref{mdg1}. In other words, we have that for all $M_j \supseteq M^*$,
\begin{align}\label{elim}
n^{1/2}\lambda^{-1/4}\left[\E(\theta|M_j,\DD)-\theta^*\right] &\xrightarrow{d} N\left(0,\mathcal{V}_{M_j}\right)
\end{align}
where $\xrightarrow{d}$ represents the convergence in distribution. Since the covariance kernel $H(t,s)$ is always positive by (C4) and $V_{M_j}$ has bounded eigenvalues for all models $M_j$ by (C3), we conclude that $\mathcal{V}_{M_j}$ is nondegenerate for all models $M_j \supseteq M^*$ and its smallest eigenvalue is bounded away from zero. Therefore \eqref{elim} implies that for any model $M_j$ in Group 1,
\begin{align*}
& \|\E(\theta|M_j,\DD)-\theta^*\| = O_p\left(n^{-1/2}\lambda^{1/4}\right), \quad
\|\E(\theta|M^*,\DD)-\theta^*\|^{-1} = O_p\left(n^{1/2}\lambda^{-1/4}\right)
\end{align*}
and hence
\begin{align}\label{mratiog1}
\frac{\|\E(\theta|M_j,\DD)-\theta^*\|}{\|\E(\theta|M^*,\DD)-\theta^*\|} & = O_p(1).
\end{align}

For any model $M_j$ in Group 2, due to the compactness of $\Theta$, we have
\begin{align}\label{mratiog2}
\frac{\|\E(\theta|M_j,\DD)-\theta^*\|}{\|\E(\theta|M^*,\DD)-\theta^*\|} & \leq 2\sup_{\theta\in \Theta}\|\theta\|\cdot O_p\left(n^{1/2}\lambda^{-1/4}\right) =  O_p\left(n^{1/2}\lambda^{-1/4}\right).
\end{align}
Note that since the total number of models is finite and does not depend on $n$, the $o_p(1)$ and $O_p(1)$ terms in all the previous expressions can always be made uniform for all models.

Now we combine \eqref{detailrate}, \eqref{ratiog1}, \eqref{ratiog2}, \eqref{mratiog1}, \eqref{mratiog2}, together with the global model selection consistency $\pi(M^*|\DD)=1+o_p(1)$ from Part (i), and obtain that
\begin{align*}
& \left\|\E(\theta|\DD)-\E(\theta|M^*,\DD)\right\|  = \left\|\sum_{M_j\neq M^*}\pi(M_j)\big[\E(\theta|M,\DD)-\E(\theta|M^*,\DD)\big]\right\| \\
&\leq \sum_{M_j\neq M^*}\pi(M_j)\cdot \big[\left\|\E(\theta|M_j,\DD)-\theta^*\right\| + \left\|\E(\theta|M^*,\DD)-\theta^*\right\|\big] \\
&\leq \sum_{M_j\in \text{Group 1}} \frac{\pi(M^*|\DD)\pi(M_j)}{\pi(M^*)}\cdot c_1\lambda^{-1/2} \cdot \left[1+O_p(1)\right]\left\|\E(\theta|M^*,\DD)-\theta^*\right\| \\
&~ + \sum_{M_j\in \text{Group 2}}  \frac{\pi(M^*|\DD)\pi(M_j)}{\pi(M^*)}\cdot c_2  e^{-\lambda\gamma/2} \cdot O_p\left(n^{1/2}\lambda^{-1/4}\right)\left\|\E(\theta|M^*,\DD)-\theta^*\right\| \\
&\leq O_p\left(\lambda^{-1/2}\right)\cdot \left\|\E(\theta|M^*,\DD)-\theta^*\right\| = o_p(1)\cdot \left\|\E(\theta|M^*,\DD)-\theta^*\right\|.
\end{align*}
The last inequality follows from the fact $n^{2/5}\prec \lambda \prec n^{2/3}$ and a comparison of the different orders of $n$.
\end{proof}

\begin{proof}[Proof of Theorem \ref{partselect}]

We use Proposition \ref{gibbs1}. In Proposition \ref{gibbs1}, we have from (A2) that the noise term $|u|\leq 2\sup_{\theta\in\Theta} |R_n(\theta)-R(\theta)|=o_p(1)$.
Now we applying (\ref{ljt}) to bound $r$. Note that $\lambda$ increases to $\infty$ due to (A1), and
$\Pi\left(\left\{\theta: R(\theta)-\inf_{\theta\in\Theta} R(\theta) <a\right\}\right)>0$ for any small $a>0$ due to (A3).
For any small $a>0$, due to (\ref{ljt}) and $\lambda \succ 1$ in (A1), we can take $\lambda > -(1/a)/ \ln \Pi\left(\left\{\theta: R(\theta)-\inf_{\theta\in\Theta} R(\theta) <a\right\}\right)$ and have that
$$ r \leq  a - \frac{1}{\lambda}\ln \Pi\left(\left\{\theta: R(\theta)-\inf_{\theta\in\Theta} R(\theta) <a\right\}\right) \leq 2a. $$
So $r=o(1)$. Together with the lower bound on the gap in (A4), we obtain from Proposition \ref{gibbs1} that
\begin{align*}
&1-\pi(M^*|\DD) = \Pi(\{M_j: j\in \Jw\} |\DD) \leq e^{-\frac{1}{2}\lambda [\gamma-o(1)-o_p(1)]} =o_p(1),
\end{align*}
since $\lambda \succ 1$ in (A1). Thus the oracle property \ref{O1} is proved. Finally we apply Proposition \ref{tvmisprob} and obtain that the oracle property \ref{O2} for BMA.
\end{proof}

\begin{proof}[Proof of Corollary \ref{marginal}]

Theorem \ref{partselect} can be applied to establish the oracle properties \ref{O1} and \ref{O2} for the combined parameter $\theta=(\phi,\omega)$. Marginalization preserves convergence in total variation norm. For the marginal distributions of $\omega$, note that
\begin{align*}
\int\left| \pi(\omega|\DD)-\pi(\omega|M^*,\DD)\right|d\omega &= \int\left| \int \pi(\omega,\phi |\DD)d\phi -  \int\pi(\omega,\phi |M^*,\DD)d\phi\right|d\omega \\
& \leq \iint\left|\pi(\omega,\phi|\DD)- \pi(\omega,\phi|M^*,\DD) \right|d\omega d\phi,
\end{align*}
which indicates that the total variation distance between the marginal densities on $\omega$ are smaller than that of the joint densities. Therefore, by Proposition \ref{tvmisprob}, the BMA oracle property also holds for the marginal posterior $\omega$, with $\int\left|\pi(\omega|\DD)-\pi(\omega|M^*,\DD)\right|d\omega=o_p(1)$.
\end{proof}

\subsection{Potential pitfalls related to partial identification}\label{pitfalls}

We will provide a simple example to show that in the context of partial identification,
the MAP model choice may not follow the oracle property \ref{O3}, and may miss the true parameter.
On the other hand, model averaging can still possess the BMA oracle property \ref{O2} and can include the true parameter.
In addition, we will show that an approach that tends to choose exclusively a ``simpler'' compatible model is not safe
with partial identification, since it can miss part of the identification region and can exclude some possible locations of the unknown true value of a point parameter.

\begin{figure}[ht]
\begin{tikzpicture}[scale=3]
\draw [line width= 1mm](3.6,0)--(4,0);
\draw [<->] (-0.1,2.2) -- (-0.1,0) -- (4.15,0);
\draw [dashed, line width=0.3mm] (0,0)--(0,0.25) -- (4,0.25)--(4,0);
\draw [line width=0.3mm] (0,0)--(3.25,0)--(3.25,2) -- (3.75,2)--(3.75,0)--(4,0);
\draw [line width=0.3mm] (3,0)--(3,2.2);
\node[below] at (0,0) {$0$};
\node[below] at (3,0) {$3$};
\node[left] at (3.6,-0.05) {$3.6$};
\node[right] at (4,-0.05) {$4$};
\node[ right] at (4.15,0) {\large$\omega$};
\node[right] at (-0.1,2.2) {Density};
\node[left] at (3,2.2) {$\pi(\omega|M_1)$};
\node[above] at (2,0.25) {$\pi(\omega|M_2)$};
\node[above] at (3.5,2) {$\pi(\omega|M_3)$};
\draw [fill](3.9,0) circle [radius=0.02];
\node[below] at (3.9,0){$\E Y$};
\end{tikzpicture}
\caption{Three models  with prior densities  $\pi(\omega|M_j)$, $j=1,2,3$.
The identification region is $[3.6,4]$ (the thick line segment). The true value of the point parameter is $\E Y$ (the dot). Model $M_1$ (with a spike prior located at $\omega=3$) is incompatible. Models $M_2$ and $M_3$ (the lower and higher plateau curves, respectively) are compatible. Of the two compatible models, Model $M_3$ has a smaller prior support but misses a
part of the identification region where the true parameter $\E Y$ is possible located.}
\label{fig.1}
\end{figure}
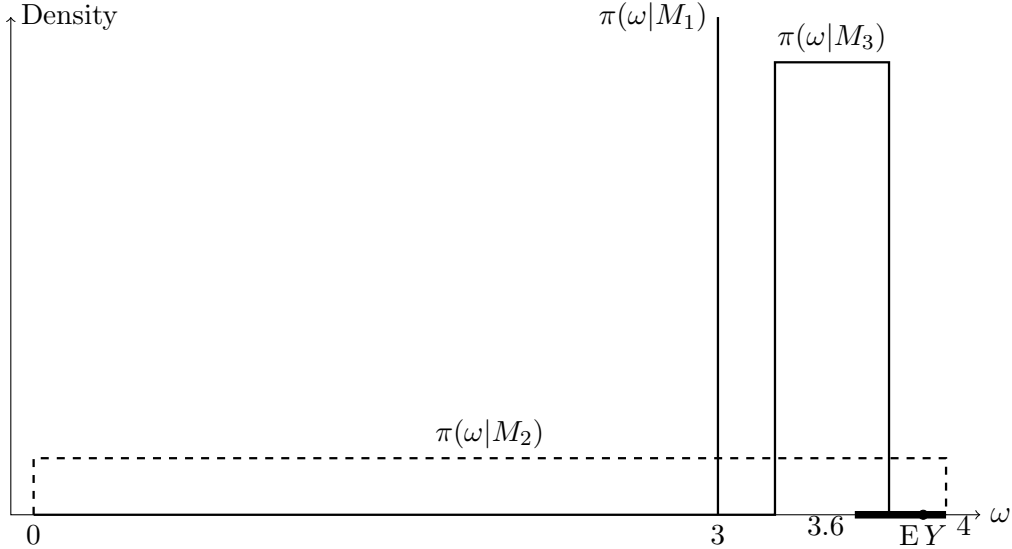

Suppose we add another ``compatible'' model in the simple example of Section \ref{simpleexample}, where $j=1,2,3$, $\pi(M_j)=1/3$ for all $j$. $\pi(\omega|M_1)=\delta_{3}(\omega)$ is a point mass supported on $W_1=\{3\}$, proposing mean GPA to be 3. $\pi(\omega|M_2)=0.25I(\{\omega\in [0,4]\})$ is a prior supported on $W_2=[0,4]$, proposing no restriction on the mean GPA. $\pi(\omega|M_3)=2 I(\{\omega\in [3.25,3.75]\})$ is a prior supported on $W_3=[3.25,3.75]$, proposing a  range of ``good'' mean GPA. The true value of $\phi$ is $\phi^*=3.6$ and the true mean GPA falls in the identification region $\Omega(\phi^*)=[3.6,4]$ as before, which intersects the supports of $\pi(\omega|M_j)$ for $j=2,3$, but not for $j=1$. Then it can be shown that models $j=2,3$ are both compatible. See Figure \ref{fig.1} for an illustration.

In this case, our oracle property \ref{O2} of Theorem \ref{partselect} can be established and applied to show that the incompatible model $M_1$ will be ignored in the limiting posterior under model averaging. Then by applying the method of Theorem 1 (ii) in \citet{ms12}, the limiting posterior of $(\omega,M_j)$ for $j=2,3$ will be $\pi(\omega,M_j|\phi=\phi^*,M_2 \text{ or } M_3) \propto \pi(\phi^*|\omega)\pi(M_j|M_2 \text{ or }M_3)\cdot$ $\pi(\omega|M_j)$. The corresponding marginal in $\omega$, $\pi(\omega |\phi=\phi^*,M_2 \text{ or } M_3)$, will be the mixture prior \\
$\sum_{j=2}^3 \pi(M_j|M_2 \text{ or }M_3) \pi(\omega|M_j)$ truncated to the identification region $[3.6,4]$, since $\pi(\phi^*|\omega)\propto I(\{\omega\in [3.6,4]\})$
for $\phi^*=3.6$ and the uniform prior $\pi(\phi|\omega)\propto I_{[\omega-1,\omega]\cap[0,4]} (\phi)$. The corresponding marginal probability of model $M_j$ ($j=2,3$) will be $\pi(M_j|\phi=\phi^*,M_2 \text{ or } M_3)\propto \pi(M_j|M_2 \text{ or } M_3) \pi(\omega\in[3.6,4] |M_j)$, i.e. $\{1/4,3/4\}$ respectively for $j=2,3$. This implies that the simpler compatible model $M_3$ is preferred, but both models $M_2$ and $M_3$ have non-vanishing posterior probabilities and will be mixed in the limit. The limiting oracle posterior of $\omega$ based on the compatible models is therefore $\pi_\infty(\omega|M_2 \text{ or } M_3) \propto 0.25I(\{\omega\in[3.6,4]\})+2I(\{\omega\in[3.6,3.75]\})$. The MAP method will tend to select $M_3$ only, since the model $M_3$ has a larger limiting posterior probability of $3/4$. The limiting posterior conditional on the MAP model choice is $\pi_\infty(\omega|M_3) \propto I(\{\omega\in[3.6,3.75]\})$. Because these limiting posteriors differ, the MAP oracle property \ref{O3} fails, but the model averaging oracle property \ref{O2} still holds.

The MAP model choice is ``simpler'' since it proposes a tighter posterior support $[3.6,3.75]$ for $\omega$. However, since it misses part of the identification region $\Omega^*=[3.6,4]$, the posterior conditional on the MAP model $M_3$ can miss some possible location of the true parameter $\omega=\E Y$ (such as 3.9). In contrast, Bayesian model averaging still observes the oracle property \ref{O2}, and its limiting posterior $\pi_\infty(\omega|M_2 \text{ or } M_3)$ has a support including the entire identification region $\Omega^*=[3.6,4]$, which does not miss any possible location of the true $\omega^*$.

In summary, in contrast to the identified models, these are two special features related to partial identification: (i) MAP model selection is no longer asymptotically equivalent to BMA, and it can miss non-unique compatible models and can miss the true parameter, while BMA is more conservative and accommodates all compatible models and all possible locations of the true parameter. (ii) More complex compatible models do not disappear asymptotically in BMA, and this is actually a good property to prevent missing possible locations of the true parameter. This property holds even when those compatible models have different dimensions, as described in the technical report \citet{jl15} Section 6.6.2. \\

\subsection{Mean oracle property with partial identification}\label{meanoracleex}

Regarding the mean oracle property \ref{O4} for the structural parameter of interest $\omega$ in the context of Section \ref{oraclepid}, we conjecture that it usually holds when we assume that $\omega$ is bounded. Since $1-\pi(M^*|\DD)=o_p(1)$  according to Theorem \ref{partselect}, we can easily apply \eqref{detailrate} to prove that $\E(\omega|\DD)-\E(\omega|M^*,\DD)=o_p(1)$. This is compared to the rate of $\E(\omega|M^*,\DD)-\omega^*$, which is typically of order 1 under partial identification. So $\E(\omega|\DD)-\E(\omega|M^*,\DD)=o_p(1)[\E(\omega|M^*,\DD)-\omega^*]$ usually holds.

We illustrate why $\E(\omega|M^*,\DD)-\omega^*$ is of order 1 in the simple example of Section \ref{simpleexample}. Models in $\Jt$ include only the full model $M_2$, and Theorem 1(ii) in \citet{ms12} can be used to show that $\int|\pi(\omega|M^*,\DD)-\pi(\omega|\phi^*,M^*)|d\omega=o_p(1)$ where $\phi^*=\E Z$. Then
\begin{align*}
&\E(\omega|M^*,\DD)-\omega^* = \int \omega \pi(\omega|M^*,\DD)d\omega-\omega^* \\
& = \int \omega\pi(\omega|\phi^*,M^*)d\omega-\omega^*  + \int \omega\left[\pi(\omega|M^*,\DD)-\pi(\omega|\phi^*,M^*)\right]d\omega \\
& \leq \int \omega\pi(\omega|\phi^*,M^*)d\omega - \omega^*  + \sup_{\omega}|\omega| \int\left|\pi(\omega|M^*,\DD)-\pi(\omega|\phi^*,M^*)\right| d\omega  \\
&=\int \omega\pi(\omega|\phi^*,M^*)d\omega - \omega^* + o_p(1),
\end{align*}
where the last equality follows from Corollary \ref{marginal} and the boundedness of $\omega$. Therefore, $\E(\omega|M^*,\DD)$ $-\omega^*= \int \omega\pi(\omega|\phi^*,M^*)d\omega - \omega^* + o_p(1)$ is of order 1, unless there is a rare coincidence that the true parameter $\omega^*$ is exactly equal to the limiting posterior mean
$\int \omega\pi(\omega|\phi^*,M^*)d\omega$. In the simple example in Section \ref{simpleexample}, $\phi^*=3.6$. Note that $\pi(\omega|\phi^*,M^*)\propto \pi(M_2)\pi(\omega|M_2)\pi(\phi^*|\omega) \propto I(\{\omega\in[\phi^*,\phi^*+1]\cap [0,4]\})$, if $\pi(\phi|\omega)\propto I(\{\phi\in [\omega-1,\omega]\})\cdot I(\{\phi\in[0,4]\})$ and $\pi(\omega|M_2)\propto I(\{\omega\in [0,4]\})$. Then the limiting posterior mean is $\int \omega\pi(\omega|\phi^*,M^*)d\omega = 3.8$. However, the true $\omega^*$ can be anywhere in $[3.6, 4]$. Unless a rare coincidence happens that $\omega^*$ is exactly $3.8$, we have that $\E(\omega|M^*,\DD)-\omega^* = 3.8-\omega^*+o_p(1)$ is of order 1, and that the mean oracle property \ref{O4} should hold.
\vspace{6mm}

\noindent {\large \bf Acknowledgement}

\vspace{3mm}

The first author is partially supported by a special fund from the Taishan Scholar Construction Project. The second author is supported by the National University of Singapore start-up grant (R155000172133).


\bibliographystyle{plainnat}
\bibliography{papers}


\end{document}